\documentclass[11pt]{article}
\usepackage{layout}
\usepackage{amssymb}
\usepackage{epsfig}
\usepackage{graphicx}
\usepackage{color}
\usepackage{latexsym}

\usepackage{amssymb}
\usepackage{algorithmicx}
\usepackage[ruled]{algorithm}
\usepackage{algpseudocode}
\usepackage{algc}

\alglanguage{pseudocode}
\usepackage{graphicx}
\usepackage{amsmath,epsfig,graphics}

\newtheorem{theorem}{Theorem}[section]

\newtheorem{remark}[theorem]{Remark}

\newcommand{\Frac}[2] {\frac{\textstyle #1} {\textstyle #2}}

\newcommand{\norm}[1]{\| #1 \|}
\newcommand{\ds}{\displaystyle}

\begin{document}

\title{Stabilized bi-grid projection methods in Finite Elements  for the 2D incompressible Navier-Stokes}
\maketitle \centerline{\scshape Hyam Abboud$^{1}$, Clara Al
Kosseifi$^{2,3}$ and Jean-Paul Chehab$^2$}
\medskip
{
  \centerline{$^{1}${\footnotesize
D\'epartement de math\'ematiques, Facult\'e des Sciences II,
Universit\'e Libanaise, Fanar, Liban}}
\centerline{$^{2}${\footnotesize Laboratoire Ami\'enois de
Math\'ematiques Fondamentales et Appliqu\'ees (LAMFA), {\small UMR
CNRS} 7352}}
  \centerline{{\footnotesize Universit\'e de Picardie Jules Verne, 33 rue Saint Leu, 80039 Amiens France} }
\centerline{$^{3}${\footnotesize Laboratoire de Physique
Appliqu\'ee (LPA), Facult\'e des Sciences II, Universit\'e
Libanaise, Fanar, Liban}}

\date{}

\begin{abstract}
\noindent
We introduce a family of bi-grid schemes in finite elements for solving 2D incompressible Navier-Stokes equations in velocity and pressure $(u,p)$.
The new schemes are based on projection methods and use two pairs of FEM spaces, a sparse and a fine one. The main computational effort is done on the coarsest velocity space with an implicit and  unconditionally  time scheme while its correction on the finer velocity space is realized with a simple  stabilized semi-implicit scheme whose the lack of stability is compensated by a high mode stabilization procedure; the pressure is updated using the free divergence property.
The new schemes are tested on the lid driven cavity up to $Re=7500$. An enhanced stability is observed as respect to classical semi-implicit methods and an important gain of CPU time is obtained as compared to implicit projection schemes.
\end{abstract}
{\small
{\bf Keywords:} {Navier-Stokes equation, bi-grid method, stabilization, Chorin-Temam projection, separation of the scales}\\
\hskip 0.2in{\bf  AMS Classification}[2010]: {35N57, \ 65L07,\ 65M60, \ 65N55}}\\
\section{Introduction}

Multigrid methods have been widely developed since more than 40 years and were proposed as fast solvers to the numerical solution of elliptic problems, see e.g. \cite{Hackbusch} for steady linear and nonlinear Dirichlet Problems in finite differences, but also for Steady Navier-Stokes Equations \cite{Bruneau1990,Ghia1982}.\\

The method is first defined for two levels of discretization (bi-grid case). The two key ingredients are both the separation of the high and the low mode components of the error provided by the use of coarse and fine grids (or spaces) $V_H$ and $V_h$ respectively,  and also the concentration of the main computational effort on the coarset (lower dimensional) subspace $V_H\subset V_h$; high modes can be represented in $V_h$ while only low modes can be captured in $V_H$. This leads to a drastic save in CPU computation time. The correction to the fine space components belonging to $V_h$ is usually realized by a simple and fast numerical scheme and is associated to a high mode smoothing.  Then the scheme is recursively applied on a set of nested grids (or spaces). \\

When considering nonlinear dissipative equations, we can distinguish two approaches for bi-grid methods:\\

In the first one,  the low and the high mode components are explicitly handled: thanks to the parabolic regularization, it is expected that the low mode components which carry the main part of the energy of the (regular) solution have a different dynamics from the high mode components, that can be seen as a correction, see \cite{DuboisJauberteauTemam} and the references therein. Also, a way to speed up the numerical integration is then to apply different schemes to these two sets of components, concentrating the effort on the computation of the low modes, that belong in $V_H$, see \cite{CDGT2001}. \\
When dealing with spectral methods, the separation in frequency is natural but it is not the case when Finite Differences or FEM methods are used for the spatial discretization. To separate the modes, a hierarchical basis approach  is used \cite{Yserentant}: using a proper interpolation (or projection) operator between $V_H$ and $V_h$, one builds a transfer operator which defines a pre-conditioner for the stiffness matrices but it allows also to express the solution in terms of main part, associated to the low mode components, and of a fluctuant part, of lower magnitude, and associated to high modes components; we refer the reader to \cite{Xu1992,Xu1994,Xu1996} and the references therein for Finite Elements discretizations, \cite{ChehabCosta2004,Pouit} for Finite Differences and \cite{BousquetMarionPetcuTemam,FaureLaminieTemam} in Finite Volumes. These schemes showed to be efficient, however they necessitate to build and manipulate hierarchical bases.\\

In the second approach, the mode separation is not used explicitly and the methods consist, at each time step, in first computing the coarse approximation $u_H$ to the fine solution $u_h$ by an unconditionally implicit stable scheme and then to update the fine space approximation by using a linearized scheme at an extrapolated value ${\tilde u}_h$ of $u_H$ in $V_h$.\\
These schemes allow to reduce the computational time with an optimal error as respected to the classical scheme when choosing accurately the mesh size of $V_h$ as respect to the one of $V_H$.
It is to be underlined that ${\tilde u}_h$ represents the mean part of the solution while $z_h=u_h-{\tilde u}_h\in V_h$ represents the fluctuant part (carrying the high mode components of the solution $u_h\in V_h$) but is not simulated in the schemes.
Such methods were developed and applied to the solution of time dependent incompressible Navier-Stokes equations
\cite{AbbGirSay2}, \cite{AbbSay}, \cite{GirLion20012} and \cite{HeLiu2006}. 
Since a linearization is used on $V_h$ the matrix to solve at each time iteration changes and it can be costly in transient regime.\\

The numerical methods we propose here are somehow in between and are inspired from the approach developed in \cite{HACKJPC} for Allen-Cahn Equation: as previously described, the use of two levels of discretization allows to concentrate the effort on the coarset, yet lower dimensional, space and at the same time to decompose the solution into its mean and fluctuant part $u_h={\tilde u}_h+z_h$. The $z_h$ are not explicitly simulated, however, they are used implicitly for a high mode stabilization as follows: consider the time integrations of the reaction-diffusion equation, in its variational form:
\begin{eqnarray}\label{reactdiff}
\Frac{d}{dt}(u,v)+(\nabla u, \nabla v)+(f(u),v)=0, \forall v \in V,\\
v(0)=v_0,
\end{eqnarray}
where $V$ is a proper Hilbert space.\\

Given two finite dimensional subspaces of $V$, $W_H$ and $V_h$, with $dim(W_H)<<dim(V_h)$, we define the bigrid-scheme as\\
\bigskip
\begin{minipage}[H]{16cm}
  \begin{algorithm}[H]
    \caption{Bi-grid Stabilized scheme for Reaction-Diffusion}\label{}
        \begin{center}
    \begin{algorithmic}[1]
        \State $u_h^{0},u^{0}_H$ given\\
            \For{$k=0,1, \cdots$}
             \State {\bf Step 1 (Coarse Space Implcit Scheme)} 
              \State {}$(\Frac{u^{k+1}_H-u^{k}_H}{\Delta t},\psi_H)+(\nabla u_H^{k+1},\nabla \psi_H)
                             +(f(u_H^{k+1}),\psi_H)=0
            , \; \forall \psi_H\in
                W_H$\label{SchemaRef}
               \State {\bf Step 2 (Fine Space semi-implicit Scheme)} 
                \State {}$(\Frac{u^{k+1}_h-u^k_h}{\Delta t},\phi_h)+{\bf \tau (u^{k+1}_h-u^k_h,\phi_h)}
               +(\nabla u^{k+1}_h,\nabla \phi_h)+
             (f(u_h^k),\phi_h)= {\bf \tau (u^{k+1}_H-u^k_H,\phi_h)}, \forall \phi_h\in V_h$
                        \EndFor
    \end{algorithmic}
     \end{center}
    \end{algorithm}
\end{minipage}
\bigskip
It is not necessary to have $W_H\subset V_h$, an inf-sup like compatibility condition has to be satisfied for defining uniquely the prolongation ${\tilde u}$, \cite{HACKJPC}.
The competition between the terms ${\bf \tau (u^{k+1}_h-u^{k}_h,\phi_h)}$ and ${\bf \tau (u^{k+1}_H-u^{k}_H,\phi_h)}$ are interpreted as a high mode filter and $\tau>0$ is the stabilization parameter  
\cite{HACKJPC} and allow to compense the explicit treatment of  the nonlinear term. Furthermore, this stabilization has few influence on the global dynamics, \cite{HACKJPC}.\bigskip

The aim of the present work is to  adapt this two-grid scheme to Navier-Stokes equations. To this end, we consider a projection method that splits the time resolution into two steps: firstly a parabolic equation on the velocity and secondly an  elliptic equation on the pressure. The bi-grid stabilization will be applied to the first step.\\
\\

The paper is organized as follows: in 
Section 2 we first present briefly the principle of the bi-grid framework (reference scheme, separation of the modes, high mode stabilization), then, after recalling the definition and some properties of the projection scheme (reference scheme), we propose different bi-grid projection schemes. Section 3 is dedicated to the numerical results. We consider the classical benchmark lid-driven cavity problem. When 
computing the steady states, the results we obtain agree with those of the literature, also an important  gain of CPU time is obtained as respect to classical methods, for a comparable precision.
The article ends in Section 4  with  concluding remarks. All the computation have been realized using FreeFem++ \cite{FreeFem}.

\bigskip

\section{Derivation of the Bi-grid Projection Schemes}
We here build the bi-grid high mode stabilized projection schemes. For that purpose, we first recall the different approaches of the bi-grid schemes in finite elements and then describe the stabilization procedure, for reaction diffusion problems. Then, we present briefly the projection schemes in finite elements we will start from to introduce  the new Bi-grid Projection Schemes.
\subsection{Principle of the bi-grid approach}
As stated in the introduction, when considering nonlinear dissipative equations, a well known way to obtain a gain in CPU time is to consider two levels of discretization and to concentrate the main computational effort on the lower 
dimensional approximation space $W_H$, while the higher accurate approximation to the solution on the fine space $V_h$ is updated by using a simple semi-implicit (yet linear) scheme.
The general pattern of a bi-grid scheme for the reaction diffusion equation (\ref{reactdiff}) writes as \\
\bigskip
\begin{minipage}[H]{16cm}
  \begin{algorithm}[H]
    \caption{Bi-grid Scheme for Reaction-Diffusion}\label{BGRD}
        \begin{center}
    \begin{algorithmic}[1]
        \State $u_h^{0},u^{0}_H$ given\\
            \For{$k=0,1, \cdots$}
             \State {\bf Step 1 (Coarse Space Implicit Scheme)} 
              \State {}$(\Frac{u^{k+1}_H-u^{k}_H}{\Delta t},\psi_H)+(\nabla u_H^{k+1},\nabla \psi_H)
                             +(f(u_H^{k+1}),\psi_H)=0
            , \; \forall \psi_H\in
                W_H$\label{SchemaRef}
               \State {\bf Step 2 (Fine Space semi-implicit Scheme)} 
                \State {}$(\Frac{u^{k+1}_h-u^k_h}{\Delta t},\phi_h)+ (\nabla u^{k+1}_h,\nabla \phi_h)+
             (f(u_h^k),\phi_h)+(f'(u_H^{k+1})u_h^{k+1},\phi_h)=0, \forall \phi_h\in V_h$
                        \EndFor
    \end{algorithmic}
     \end{center}
    \end{algorithm}
\end{minipage}
\\
This approach was proposed by \cite{GirLion20012,HeLiu2006,Layton1993}. However, $f'(u_H^{k+1})$, the linearized of $f$ at $u_H^{k+1}$ must be computed at each time step, changing the matrix to solve at each iteration. A way to avoid this drawback is to consider a classical semi-implicit scheme and to compense the lack of stability by adding, e.g., a first order damping term $\tau (u^{k+1}_h-u^k_h)$, obtaining as second step:\\
\begin{minipage}[H]{16cm}
  \begin{algorithm}[H]
      \begin{algorithmic}[1]
\State {\bf Step 2 (Stabilized Fine Space semi-implicit Scheme)} 
                \State {}$(\Frac{u^{k+1}_h-u^k_h}{\Delta t},\phi_h)+ \tau (u^{k+1}_h-u^k_h,\phi_h)       +(\nabla u^{k+1}_h,\nabla \phi_h)+
             (f(u_h^k),\phi_h)=0, \forall \phi_h\in V_h$
    \end{algorithmic}
    \end{algorithm}
    \end{minipage}
    \\
    This  stabilization procedure allows to obtain unconditionally stable time scheme for large values of $\tau>0$, that can be tuned. The resulting scheme is fast, however, it can slow down drastically the dynamics, particularly convergence to steady states occurs in longer times. This is due to the fact that the damping acts on all the mode components, including the low ones which are associated to the mean part of the function and carry the main part of the $L^2$ energy, when considering a Fourier-like interpretation. A way to overcome this drawback is to apply the stabilization to the only high mode components 
 $z_h$ of $u_h$ which correspond to a fluctuent part of $u_h$ ; the stability of a scheme is indeed related to its capability to contain the propagation of the high frequencies.  Hence the scheme becomes\\
 \\
 \begin{minipage}[H]{16cm}
  \begin{algorithm}[H]
      \begin{algorithmic}[1]
\State {\bf Step 2 (High Mode Stabilized Fine Space semi-implicit Scheme)} 
                \State {}$(\Frac{u^{k+1}_h-u^k_h}{\Delta t},\phi_h)+ \tau ({\cal P}( u^{k+1}_h-u^k_h),\phi_h)       +(\nabla u^{k+1}_h,\nabla \phi_h)+
             (f(u_h^k),\phi_h)=0, \forall \phi_h\in V_h$
    \end{algorithmic}
    \end{algorithm}
    \end{minipage}
    \\
    \\
 where ${\cal P}( u^{k+1}_h-u^k_h)$ capture the high mode components of $u^{k+1}_h-u^k_h$.\\
 
 At this point we can distinguish two main ways to decompose $u_h$ as a sum of its mean part ${\tilde u}_h$ and its fluctuent part $z_h$ in finite elements:
\begin{itemize}
\item Hierarchical basis in Finite Elements: the fine space FEM approximation space $V_h$ is decomposed as $V_h=V_{H}\bigoplus W_h\subset H^1$, where $V_H\subset V_h$ is the coarse FEM space that can capture only low modes, and $W_h$ is the complementary space, generated by the basis functions of $V_h$ that do not belong in $V_H$,  and that capture high mode components. When using $\mathbb{P}_1$ finite elements, the components of $W_h$ of a function $u_h\in V_h$  are built as proper local interpolation error from $u_H\in V_H$. The functions of approximation $V_h$ can be uniquely decomposed as
$$
u_h=u_H+z_h,
$$
with $u_H\in V_H$ and $z_h\in W_h$. It has been showed that the linear change of variable $S: u_h\rightarrow (u_H,z_h)$ provides an efficient preconditioner for the stiffness matrices but also proceeds to a scale separation \cite{CalgaroLaminieTemam,CalgaroDebusscheLaminie,MarionTemam1,MarionTemam2,Yserentant}. 
\\
The spatial discretization reads then to a differential system satisfied by $u_H$ and $z_h$, namely
\begin{eqnarray}
\Frac{d}{dt}(u_H+z_h,\phi_H)+(\nabla (u_H+z_h), \nabla \phi_H)+(f(u_H+z_h),\phi_H)=0, \forall \phi_H \in V_H,\\
\Frac{d}{dt}(u_H+z_h,\psi_h)+(\nabla (u_H+z_h), \nabla \psi_h)+(f(u_H+z_h),\psi_h)=0, \forall \psi_h \in V_h.
\end{eqnarray}

New time marching methods are obtained applying two different schemes for $u_h$ and for
$z_h$; particularly, the separation of the scales allows to use a simple and fast scheme for the $z_h$ components and to save computational time, with a good accuracy. Approximations and simplifications can be applied to each equation, particularly using an approached law linking high and low mode components of type $\Phi(u_H)=z_h$, \cite{MarionTemam1,MarionTemam2}.\\

The approach followed in \cite{MarionXu} is technically different but is based on the same principle: Denoting by $P_H$ the $L^2$ orthogonal projection from $V_h$ on $V_H$, the decomposition
$u_h=u_H+w_h$ with $u_H=P_Hu_h$ and $w_h=(Id-P_H)u_h$ is used, making $u_H$ and $w_h$ orthogonal. They proposed and analyzed schemes of Nonlinear Galerkin type, says in which 
an asymptotic approach law $\Phi(u_H)=w_h)$ is implemented together with a time marching scheme on $u_H$.
A drawback is that one must build a basis for $(Id-P_H)V_h$.
\\
Finally let us mention similar developments in finite differences \cite{CCLZ,ChehabCosta2004} and spectral methods \cite{CDGT2001,DuboisJauberteauTemam,DuboisJauberteauTemamTribbia,JauberteauTemamTribbia}
\item A $L^2$-like filtering, used in \cite{HACKJPC} and to which we will concentrate: we define first the prolongation operator ${\cal P}:
V_H\rightarrow V_h$ by
\begin{eqnarray}\label{uprolong}
(u_H-{\cal P} (u_H), \phi_h)=0, \forall \phi_h \in V_h.
\end{eqnarray}
Then, setting ${\tilde u}_h={\cal P} (u_H)$, we write
$$
u_h={\tilde u}_h+(u_h-{\tilde u}_h)={\tilde u}_h+z_h.
$$
\end{itemize}
We now apply the high mode stabilization to the $z_h$ components and get, after usual simplifications\\
\\
\begin{minipage}[H]{16cm}
  \begin{algorithm}[H]
      \begin{algorithmic}[1]
\State {\bf Step 2 (High Mode Stabilized Fine Space semi-implicit Scheme)} 
                \State {}$(\Frac{u^{k+1}_h-u^k_h}{\Delta t},\phi_h)+ \tau ( z^{k+1}_h-z^k_h,\phi_h)       +(\nabla u^{k+1}_h,\nabla \phi_h)+
             (f(u_h^k),\phi_h)=0, \forall \phi_h\in V_h$                     
    \end{algorithmic}
    \end{algorithm}
    \end{minipage}
    \\
and using the identity 
$$
(z^{k+1}_h-z^k_h,\phi_h)=(u^{k+1}_h-u^k_h,\phi_h)-({\tilde u}^{k+1}_h-{\tilde u}^k_h,\phi_h)=(u^{k+1}_h-u^k_h,\phi_h)-(u^{k+1}_H-u^k_H,\phi_h), \forall \phi_h\in V_h
$$
we can write the high mode stabilized bi-grid scheme as\\
\bigskip
\begin{minipage}[H]{16cm}
  \begin{algorithm}[H]
    \caption{Bi-grid Stabilized scheme for Reaction-Diffusion}\label{BGSTRD}
        \begin{center}
    \begin{algorithmic}[1]
        \State $u_h^{0},u^{0}_H$ given\\
            \For{$k=0,1, \cdots$}
             \State {\bf Step 1 (Coarse Space Implicit Scheme)} 
              \State {}$(\Frac{u^{k+1}_H-u^{k}_H}{\Delta t},\psi_H)+(\nabla u_H^{k+1},\nabla \psi_H)
                             +(f(u_H^{k+1}),\psi_H)=0
            , \; \forall \psi_H\in
                W_H$\label{SchemaRef}
               \State {\bf Step 2 (Fine Space semi-implicit Scheme)} 
                \State {}$(\Frac{u^{k+1}_h-u^k_h}{\Delta t},\phi_h)+{\bf \tau (u^{k+1}_h-u^k_h,\phi_h)}
               +(\nabla u^{k+1}_h,\nabla \phi_h)+
             (f(u_h^k),\phi_h)= {\bf \tau (u^{k+1}_H-u^k_H,\phi_h)}, \forall \phi_h\in V_h$
                        \EndFor
    \end{algorithmic}
     \end{center}
    \end{algorithm}
\end{minipage}
\bigskip

We now consider incompressible Navier-Stokes equations. FEM Bi-grid methods for incompressible NSE as proposed, e.g.,  in  \cite{AbbGirSay2,AbbSay,GirLion20012,HeLiu2006,Layton1993}, are built using the pattern of scheme \ref{BGRD}, and written as following. Given two pairs of FEM spaces $(X_H,Y_H)$ and $(X_h, Y_h)$ satisfying the discrete inf-sup condition, with $X_H\subset X_h$ and $Y_H\subset Y_h$, one defines  the bi-grid iterations as
\bigskip
\begin{minipage}[H]{16cm}
  \begin{algorithm}[H]
    \caption{Bi-grid Scheme for Navier-Stokes Equation}\label{BGNS}
        \begin{center}
    \begin{algorithmic}[1]
        \State $u_h^{0},u^{0}_H$ given\\
            \For{$k=0,1, \cdots$}
             \State {\bf Step 1 (Coarse Space Implicit Scheme)} 
              \State {}$(\Frac{u^{k+1}_H-u^{k}_H}{\Delta t},\psi_H)+\Frac{1}{Re}(\nabla u_H^{k+1},\nabla \psi_H)
                             +((u_H^{k+1}.\nabla) u_H^{k+1},\psi_H)-(div (\psi_H),p^{k+1}_H)=(f,\psi_H), \forall \psi_H\in X_H$
               \State $(div (u^{k+1}_H),q_H)=0, \forall q_H\in Y_H$
                \label{SchemaRef}
               \State {\bf Step 2 (Fine Space semi-implicit Scheme)} 
                \State {}$(\Frac{u^{k+1}_h-u^k_h}{\Delta t},\phi_h)+ \Frac{1}{Re} (\nabla u^{k+1}_h,\nabla \phi_h)+
          (u_H^{k+1}.\nabla u^{k+1}_h,\phi_h)-(div (\psi_h),p^{k+1}_h)=(f,\psi_h), \forall \psi_h\in X_h$
           \State $(div (u^{k+1}_h),q_h)=0, \forall q_h\in Y_h$

                        \EndFor
    \end{algorithmic}
     \end{center}
    \end{algorithm}
\end{minipage}
We have used here the simple linearization of the nonlinear term $(u_H^{k+1}.\nabla u^{k+1}_h,\phi_h)$ proposed in \cite{GirLion20012}, but other linearisations can be considered. At each step the method needs to solve first a mixed nonlinear FEM problem on the coarse pair FEM spaces  $(X_H,Y_H)$, then a linear mixed FEM problem on the fine FEM spaces $(X_h,Y_h)$; in this last step, the underlying matrix
changes at each iteration. To avoid this drawback and to apply the stabilized bi-grid method, as described in Algorithm \ref{BGSTRD}, we will use a projection method. In such a case the computation of the velocity is decoupled from the one of the pressure, the intermediary velocity satisfies a nonlinear convection-diffusion equation. The main idea is then to apply a stabilized bi-grid scheme to this equation to speed up the resolution. We present below several options of this approach.%
\subsection{Bi-grid Projection Schemes in Finite Elements}
First of all, let us recall the framework of the Projection Schemes in Finite elements and then derive the new bi-grid stabilized methods.\\

The mixed variational formulation of the motion
of a viscous and incompressible fluid in a domain $\Omega$ is described by
the unsteady Navier-Stokes equation
\[\left\{
\begin{array}{ll}
\medskip
(\ds\frac{\partial u}{\partial t},v)+(u\cdot \nabla u,v)-\nu (\nabla u,\nabla v)-(p,\mbox{div}\;v)=(f,v),\;\forall v\in X=(H_0^1(\Omega))^d \\
\medskip
(\mbox{div}\;u,q)=0,\;\forall q\in Y=L^2(\Omega) \\
\end{array}
\right.\] where $\Omega$ is a domain of $\mathbb{R}^d, d=2, 3$ of
lipschizian bound $\partial \Omega, \overrightarrow{n}$ the normal
outside unit and a time interval $[0,T],T>0$.\\
\noindent Here is a choice of the approximation spaces $X_h$ and
$Y_h$ for the velocity and the pressure
\begin{equation}\label{EFP2NS}
X_h=\{v\in \mathcal{C}^0(\bar{\Omega})^2|
v_{|\kappa}\in\mathbb{P}_2^2,\;\forall \kappa\in T_h,\; v=0\text{ on
}\partial\Omega\},
\end{equation}
\begin{equation}\label{EFP1NS}
Y_h=\{q\in \mathcal{C}^0(\bar{\Omega})|
q_{|\kappa}\in\mathbb{P}_1,\;\forall \kappa\in T_h\}\cap L^2_0(\Omega).
\end{equation}
%
%
\noindent The degrees of freedom for the velocity are the vertices of
the triangulation and the midpoints of the edges of the triangles 
of the triangulation $T_h$. The degrees of freedom for the pressure are the vertices
of $T_h$ assumed to be uniformly regular. It is
well known that because of the constraint of incompressibility,
the choice of $X_h$ and $Y_h$ is not arbitrary. They must satisfy
a suitable compatibility condition, the "inf-sup" condition of
Babu\u{s}ka-Brezzi (cf. \cite{Babuska1973,Brezzi1974}):
\begin{equation}\label{InfSupdis}
\ds\inf_{q\in Y_h}\sup_{u\in X_h}\ds\frac{\int_\Omega q(div\;
u)dx}{\norm{q}_{L^2(\Omega)}\norm{\nabla
u}_{L^2(\Omega)}}\geq\beta^\ast,
\end{equation}
where $\beta^\ast$ is independent of $h$. All the results
presented in this paper are done with the Taylor-Hood finite
element $\mathbb{P}_2/\mathbb{P}_1$.\\

Let us consider the semi-discretization in time and focus on
marching schemes. Let $u^k\simeq u(x,k \Delta t)$ be a sequence of
functions; $\Delta t$ is the time step. We compare our method to
the projection method applied on the fine grid. It is a fractional
step-by-step method of decoupling the computation of the velocity
from that of the pressure, by first solving a convection-diffusion
problem such that the resulting velocity is not necessarily zero
divergence; then in a second step, the latter is projected onto a
space of functions with zero divergence in order to satisfy the
incompressibility condition (cf.
\cite{Chorin,GuermondShen,GuerShen2003,GuerShen2001,Tem84,Tem6970,Tem68}).
We start by presenting the implicit reference scheme
\cite{Tem77}:
\begin{center}
\begin{minipage}[H]{16cm}
  \begin{algorithm}[H]
    \caption{Reference Scheme}\label{ALGO1G1}
    \begin{algorithmic}[1]
     \For{$k=0,1, \cdots$}
           \State {\bf Find $u_h^\ast$ in $X_h$}
           ($\ds\frac{u_h^\ast-u_h^k}{\Delta t},\psi_h) + \nu(\nabla
u_h^\ast,\nabla \psi_h)+((u_h^\ast\cdot\nabla)u_h^\ast,\psi_h)
=(f,\psi_h), \forall \psi_h \in X_h$
\State {\bf Find $p_h^{k+1}$ in $Y_h$} $(div\;\nabla
p_h^{k+1},\chi_h)= (\ds\frac{div\;u_h^\ast}{\Delta t},\chi_h),
\forall \chi_h \in Y_h$
\State {\bf Find $u_h^{k+1}$ in $X_h$} ($u_h^{k+1}-u_h^\ast+\Delta
t \nabla p_h^{k+1},\psi_h)=0, \forall \psi_h \in X_h$
            \EndFor
 \end{algorithmic}
   \end{algorithm}
\end{minipage}
\end{center}
This will be our reference scheme, used on the coarse FEM Space. It enjoys of unconditional stability properties, see \cite{Tem77}, chapter III, section 7.3.
\begin{remark}
We could use also an incremental projection method introduced by Goda \cite{Goda1979}, see also 
\cite{GuermondShen}. Goda has proven that adding a previous value
of the gradient of the pressure ($\nabla p^k$) in the first step
of the projection method then rectify the value of the velocity in
the second step will improve the accuracy, in other words, he
proposed the following algorithm:
\begin{center}
\begin{minipage}[H]{16cm}
  \begin{algorithm}[H]
    \caption{Incremental Reference}\label{ALGO1G2}
    \begin{algorithmic}[1]
     \For{$k=0,1, \cdots$}
           \State {\bf Find $u_h^\ast$ in $X_h$}  \\
           $(\ds\frac{u_h^\ast-u_h^k}{\Delta t},\psi_h) + \nu(\nabla
u_h^\ast,\nabla \psi_h)+((u_h^\ast\cdot\nabla)u_h^\ast,\psi_h)+(\nabla
p_h^k,\psi_h) =(f,\psi_h), \forall \psi_h \in X_h$
\State {\bf Find $p_h^{k+1}$ in $Y_h$} $\alpha(div\;\nabla
p_h^{k+1}-div\;\nabla p_h^k,\chi_h)=
(\ds\frac{div\;u_h^\ast}{\Delta t},\chi_h), \forall \chi_h \in
Y_h$
\State {\bf Find $u_h^{k+1}$ in $X_h$}
($u_h^{k+1}-u_h^\ast+\alpha\Delta t (\nabla p_h^{k+1}-\nabla
p_h^k),\psi_h)=0, \forall \psi_h \in X_h$
            \EndFor
 \end{algorithmic}
   \end{algorithm}
\end{minipage}
\end{center}
The results obtained are similar to the one produced by the non incremental scheme to which we focus for a sake of simplicity.
\end{remark}

We can now derive the bi-grid schemes.

\noindent Our first approach consists in separating in scales the
intermediate velocity $u^\ast_h$ which we introduce according to
the projection method of Chorin-Temam
\cite{Chorin,Tem84,Tem6970,Tem68}. 
First we compute $u_H^\ast$ on the coarse space $X_H$ and then
stabilize the high frequencies of $u_h^\ast$ on the fine space
$X_h$. Based on the free divergence condition, we find the
pressure $p_h$ whose mean is equal to zero. We end up finding $u_h$ and
restricting it to the coarse grid to get $u_H^k$ for the next
iteration.
\begin{center}
\begin{minipage}[H]{16cm}
  \begin{algorithm}[H]
    \caption{Two grids Algo1}\label{NSALGO1}
    \begin{algorithmic}[1]
     \For{$k=0,1, \cdots$}
           \State {\bf Find $u_H^\ast$ in $X_H$} ($\ds\frac{u_H^\ast-u_H^k}{\Delta t},\psi_H) + \nu(\nabla
u_H^\ast,\nabla \psi_H)+((u_H^\ast\cdot\nabla)u_H^\ast,\psi_H)
=(f,\psi_H), \forall \psi_H \in X_H$
\State {\bf Find $u_h^\ast$ in $X_h$}
            $$
(1+\tau \Delta t)(\ds\frac{u_h^\ast-u_h^k}{\Delta t},\phi_h)
+\nu (\nabla u_h^\ast,\nabla \phi_h)
+((u_h^k\cdot\nabla)u_h^k,\phi_h)
=\tau(u_H^\ast-u_H^k,\phi_h)+(f,\phi_h),\forall
\phi_h \in X_h
$$
\State {\bf Find $p_h^{k+1}$ in $Y_h$} $(div\;\nabla
p_h^{k+1},\chi_h)= (\ds\frac{div\;u_h^\ast}{\Delta t},\chi_h),
\forall \chi_h \in Y_h$
\State {\bf Find $u_h^{k+1}$ in $X_h$} ($u_h^{k+1}-u_h^\ast+\Delta
t \nabla p_h^{k+1},\phi_h)=0, \forall \phi_h \in X_h$
\State {\bf Solve in $X_H$} $(u_H^{k+1}-u_h^{k+1},\psi_H)=0,
\forall \psi_H \in X_H$
            \EndFor
    \end{algorithmic}
    \end{algorithm}
\end{minipage}
\end{center}
%
We propose in the second algorithm to replace
$(u_h^k.\nabla)u_h^k$ by $(u_H^\ast\cdot\nabla)u_h^\ast$. We obtain:
%
\begin{center}
\begin{minipage}[H]{16cm}
  \begin{algorithm}[H]
    \caption{Variant of the Two grids Algo1~: Two grids Algo2}\label{NSALGO2}
    \begin{algorithmic}[1]
     \For{$k=0,1, \cdots$}
           \State {\bf Find $u_H^\ast$ in $X_H$} ($\ds\frac{u_H^\ast-u_H^k}{\Delta t},\psi_H) + \nu(\nabla
u_H^\ast,\nabla \psi_H)+((u_H^\ast\cdot\nabla)u_H^\ast,\psi_H)
=(f,\psi_H), \forall \psi_H \in X_H$

\State {\bf Find $u_h^\ast$ in $X_h$ }
             $$
(1+\tau \Delta t)(\ds\frac{u_h^\ast-u_h^k}{\Delta t},\phi_h)
+\nu(\nabla u_h^\ast,\nabla \phi_h)
+((u_H^\ast\cdot\nabla)u_h^\ast,\phi_h)
=\tau(u_H^\ast-u_H^k,\phi_h)+(f,\phi_h),\forall
\phi_h \in X_h
$$
\State {\bf Find $p_h^{k+1}$ in $Y_h$ } $(div\;\nabla
p_h^{k+1},\chi_h)= (\ds\frac{div\;u_h^\ast}{\Delta t},\chi_h),
\forall \chi_h \in Y_h$
\State {\bf Find $u_h^{k+1}$ in $X_h$ }($u_h^{k+1}-u_h^\ast+\Delta
t \nabla p_h^{k+1},\phi_h)=0, \forall \phi_h \in X_h$
\State {\bf Solve in $X_H$} $(u_H^{k+1}-u_h^{k+1},\psi_H)=0,
\forall \psi_H \in X_H$
            \EndFor
    \end{algorithmic}
    \end{algorithm}
\end{minipage}
\end{center}
The second approach that we consider in the following is to find directly
$u^{k+1}_H$ by the projection method applied on the coarse grid. By this
technique, we do not need to restrict $u_h^{k+1}$ to $X_H$ to define $u_H^{k+1}$.
\begin{center}
\begin{minipage}[H]{16cm}
  \begin{algorithm}[H]
    \caption{Two grids Algo3 }\label{NSALGO3}
    \begin{algorithmic}[1]
     \For{$k=0,1, \cdots$}
           \State {\bf Find $u_H^\ast$ in $X_H$ }($\ds\frac{u_H^\ast-u_H^k}{\Delta t},\psi_H) + \nu(\nabla
u_H^\ast,\nabla \psi_H)+((u_H^\ast\cdot\nabla)u_H^\ast,\psi_H)
=(f,\psi_H), \forall \psi_H \in X_H$
\State {\bf Find $p_H^{k+1}$ in $Y_H$ }$(div\;\nabla
p_H^{k+1},\chi_H)= (\ds\frac{div\;u_H^\ast}{\Delta t},\chi_H),
\forall \chi_H \in Y_H$
\State {\bf Find $u_H^{k+1}$ in $X_H$ }($u_H^{k+1}-u_H^\ast+\Delta
t \nabla p_H^{k+1},\psi_H)=0, \forall \psi_H \in X_H$
\State {\bf Find $u_h^\ast$ in $X_h$ }
$$
(1+\tau \Delta t)(\ds\frac{u_h^\ast-u_h^k}{\Delta t},\phi_h)
+\nu(\nabla u_h^\ast,\nabla \phi_h)
+((u_h^\ast\cdot\nabla)u_h^\ast,\phi_h)
=\tau(u_H^\ast-u_H^k,\phi_h)+(f,\phi_h),\forall
\phi_h \in X_h
$$
\State {\bf Find $p_h^{k+1}$ in $Y_h$ }$(div\;\nabla
p_h^{k+1},\chi_h)= (\ds\frac{div\;u_h^\ast}{\Delta t},\chi_h),
\forall \chi_h \in Y_h$
\State {\bf Find $u_h^{k+1}$ in $X_h$ }($u_h^{k+1}-u_h^\ast+\Delta
t \nabla p_h^{k+1},\phi_h)=0, \forall \phi_h \in X_h$
            \EndFor
               \end{algorithmic}
    \end{algorithm}
\end{minipage}
\end{center}
%
\begin{center}
\begin{minipage}[H]{16cm}
  \begin{algorithm}[H]
    \caption{ Variant of the Two grids Algo3~: Two grids Algo4 }\label{NSALGO4}
    \begin{algorithmic}[1]
     \For{$k=0,1, \cdots$}
           \State {\bf Find $u_H^\ast$ in $X_H$} ($\ds\frac{u_H^\ast-u_H^k}{\Delta t},\psi_H) + \nu(\nabla
u_H^\ast,\nabla \psi_H)+((u_H^\ast\cdot\nabla)u_H^\ast,\psi_H)
=(f,\psi_H), \forall \psi_H \in X_H$
\State {\bf Find $p_H^{k+1}$ in $Y_H$} $(div\;\nabla
p_H^{k+1},\chi_H)= (\ds\frac{div\;u_H^\ast}{\Delta t},\chi_H),
\forall \chi_H \in Y_H$
\State {\bf Find $u_H^{k+1}$ in $X_H$} ($u_H^{k+1}-u_H^\ast+\Delta
t \nabla p_H^{k+1},\psi_H)=0, \forall \psi_H \in X_H$
\State {\bf Find $u_h^\ast$ in $X_h$}
             $$
(1+\tau \Delta t)(\ds\frac{u_h^\ast-u_h^k}{\Delta t},\phi_h)
+\nu(\nabla u_h^\ast,\nabla \phi_h)
+((u_H^\ast\cdot\nabla)u_h^\ast,\phi_h) =\tau(u_H^\ast-u_H^k,\phi_h)+(f,\phi_h),\forall
\phi_h \in X_h
$$
\State {\bf Find $p_h^{k+1}$ in $Y_h$} $(div\;\nabla
p_h^{k+1},\chi_h)= (\ds\frac{div\;u_h^\ast}{\Delta t},\chi_h),
\forall \chi_h \in Y_h$
\State {\bf Find $u_h^{k+1}$ in $X_h$} ($u_h^{k+1}-u_h^\ast+\Delta
t \nabla p_h^{k+1},\phi_h)=0, \forall \phi_h \in X_h$
            \EndFor
    \end{algorithmic}
    \end{algorithm}
\end{minipage}
\end{center}
\section{Numerical results}
\subsection{The Lid Driven Cavity}
\noindent We simulate the Navier-Stokes equations in 2D by a
stabilization technique. We use the unit square mesh
$[0,1]\times[0,1]$ and vary the dimensions of the coarse and fine FEM spaces. We choose as spaces of approximation the well-known
Taylor-Hood element $\mathbb{P}_2/\mathbb{P}_1$ for the velocity
and the pressure respectively. Here are the results of the driven
cavity flow obtained by adopting the variational formulation
provided by the projection technique (Chorin-Temam) proposed by
the two preceeding steps. In this case the velocity is imposed only
in the upper boundary with $u=(1,0)$ and zero Dirichlet conditions are
imposed on the rest of the boundary (see Fig. \ref{DrivenCavity}), below
\begin{center}
\begin{figure}[!h]
\vspace{-.5cm}
\begin{center}
\includegraphics[height=4.5cm,width=8cm]{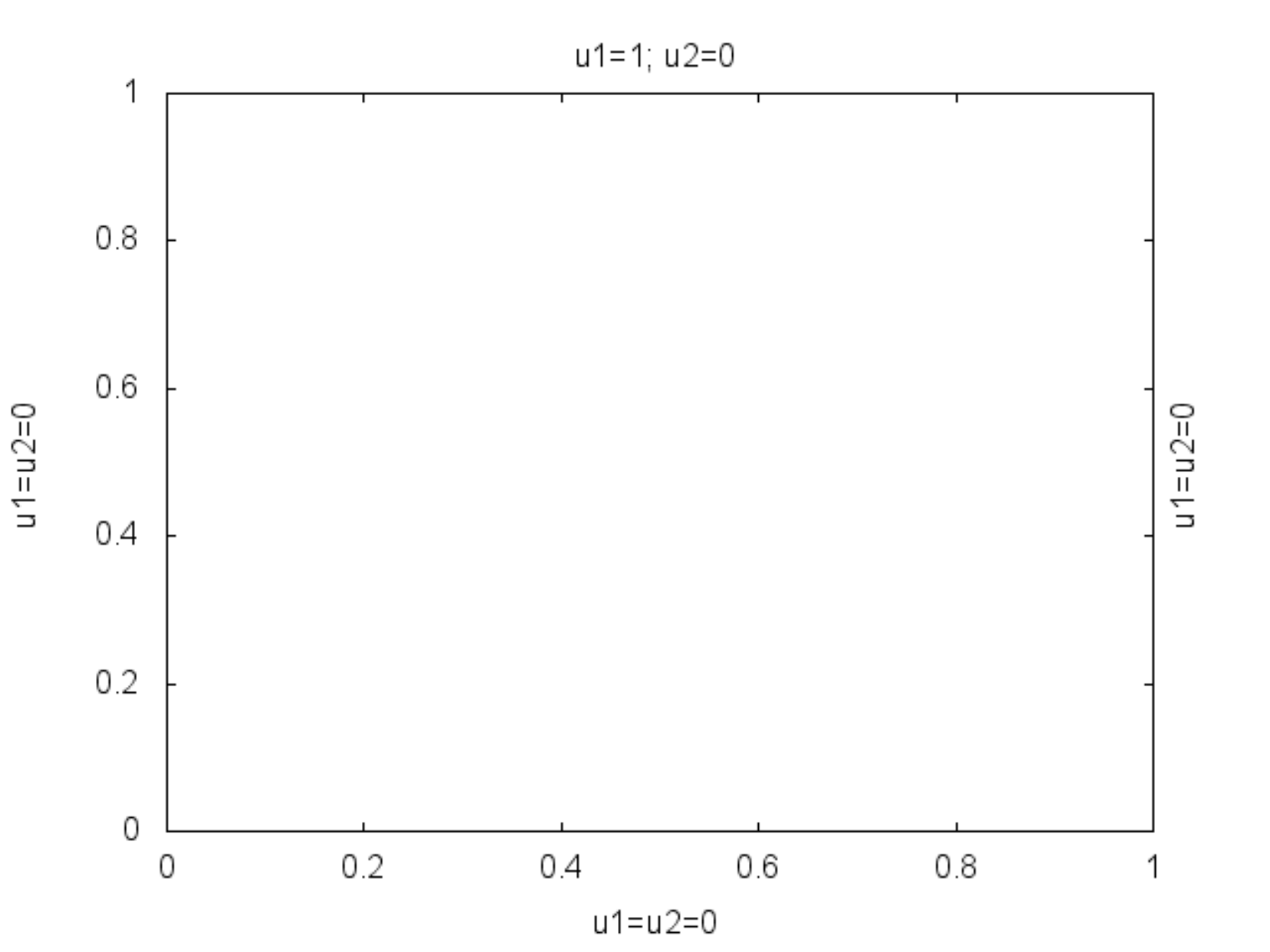}
\end{center}
\vspace{-0.3cm}
\caption{Boundary conditions}\label{DrivenCavity}
\end{figure}
\end{center}
\subsubsection{Computation of Steady States}
\begin{center}
\vspace{-0.5cm}
\begin{figure}[!h]
\begin{center}
\hspace{-0.4cm}
\includegraphics[height=3.5cm,width=4.2cm]{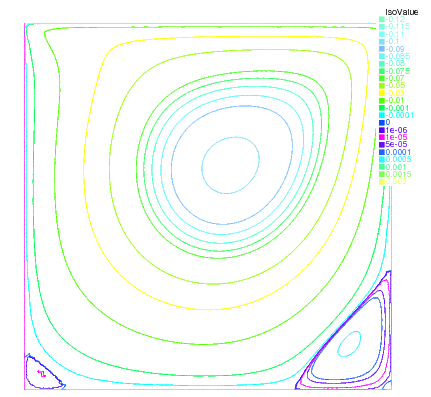}
\hspace{-0.5cm}
\includegraphics[height=3.5cm,width=4.2cm]{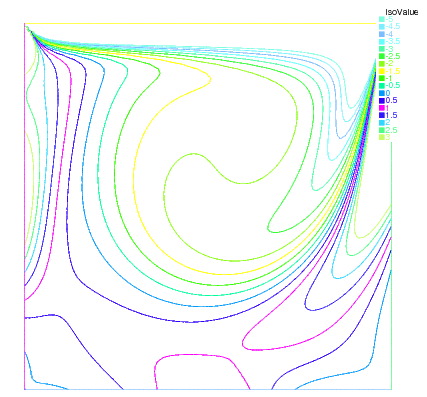}
\hspace{-0.5cm}
\includegraphics[height=3.5cm,width=4.2cm]{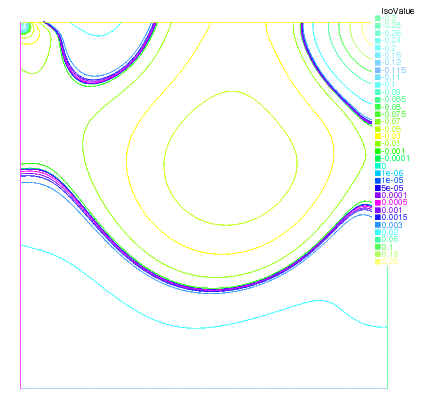}
\hspace{-0.5cm}
\includegraphics[height=3.5cm,width=4.2cm]{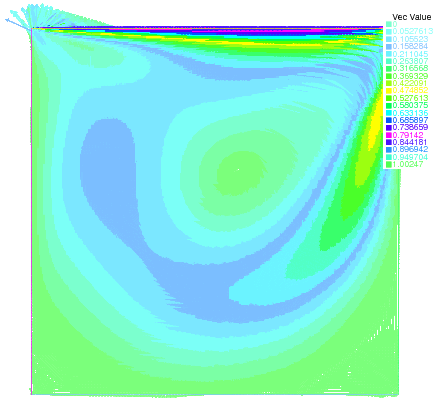}\\
\hspace{-0.4cm}
\includegraphics[height=3.5cm,width=4.2cm]{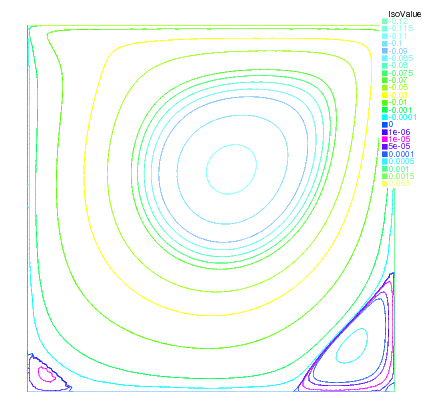}
\hspace{-0.5cm}
\includegraphics[height=3.5cm,width=4.2cm]{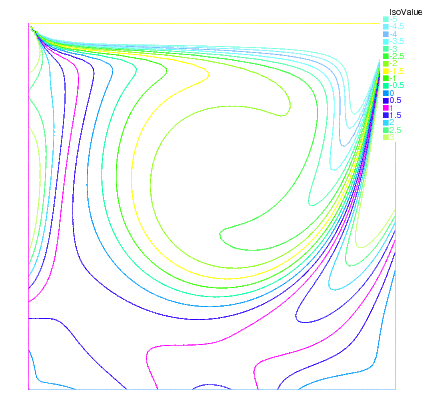}
\hspace{-0.5cm}
\includegraphics[height=3.5cm,width=4.2cm]{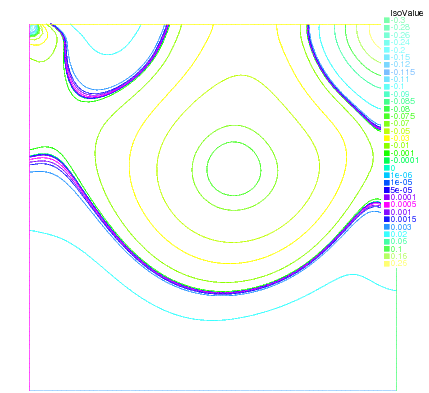}
\hspace{-0.5cm}
\includegraphics[height=3.5cm,width=4.2cm]{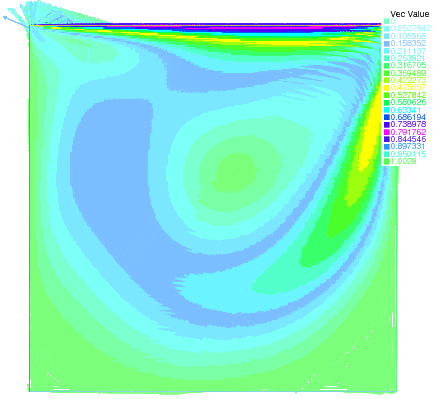}\\
\hspace{-0.4cm}
\includegraphics[height=3.5cm,width=4.2cm]{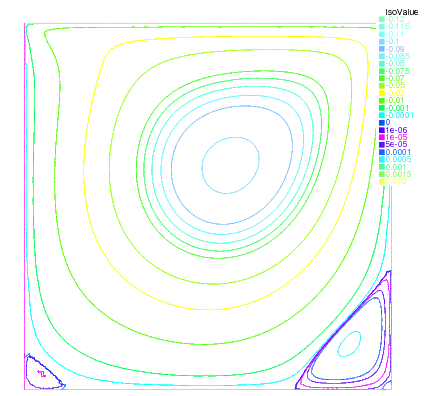}
\hspace{-0.5cm}
\includegraphics[height=3.5cm,width=4.2cm]{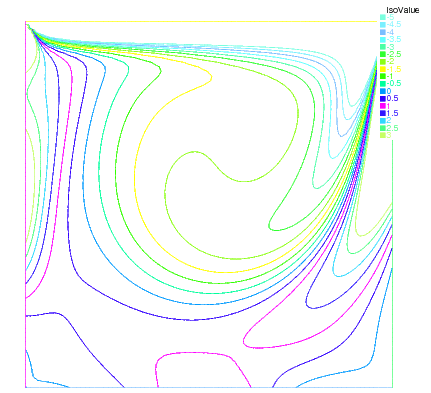}
\hspace{-0.5cm}
\includegraphics[height=3.5cm,width=4.2cm]{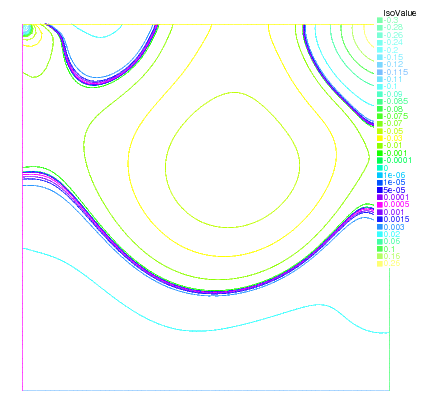}
\hspace{-0.5cm}
\includegraphics[height=3.5cm,width=4.2cm]{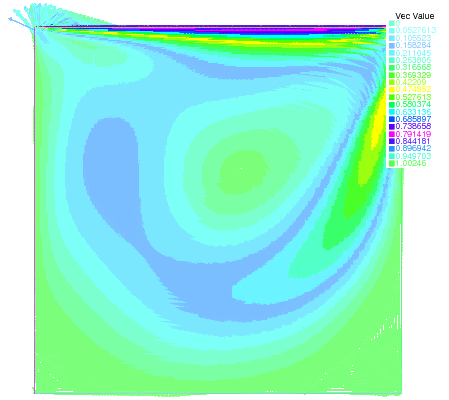}\\
\hspace{-0.4cm}
\includegraphics[height=3.5cm,width=4.2cm]{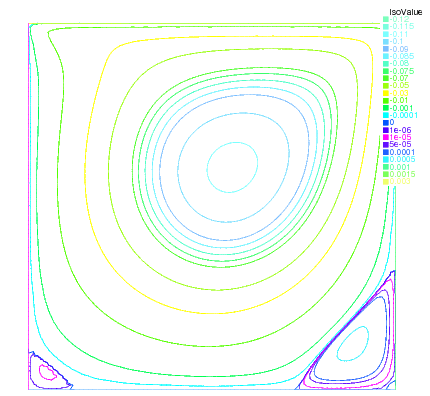}
\hspace{-0.5cm}
\includegraphics[height=3.5cm,width=4.2cm]{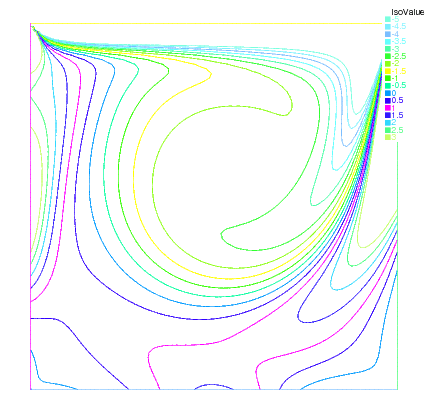}
\hspace{-0.5cm}
\includegraphics[height=3.5cm,width=4.2cm]{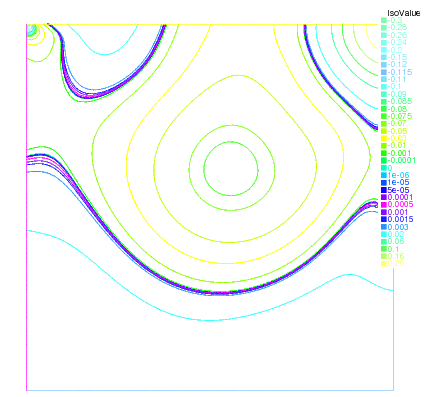}
\hspace{-0.5cm}
\includegraphics[height=3.5cm,width=4.2cm]{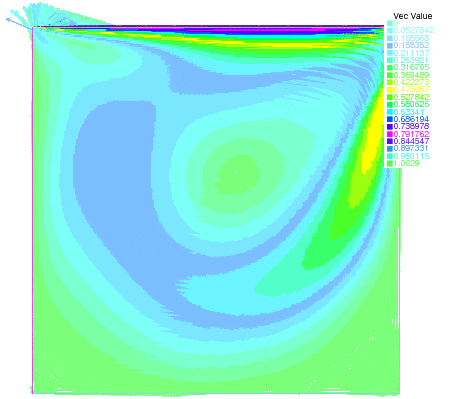}\\
%
\hspace{-0.4cm}
\includegraphics[height=3.5cm,width=4.2cm]{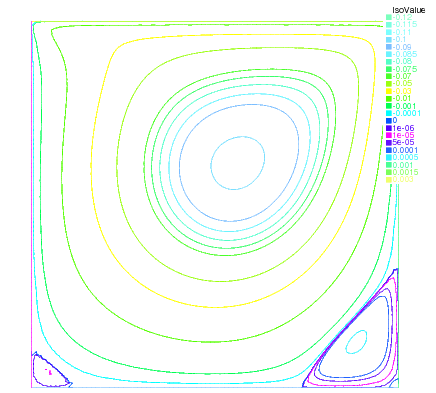}
\hspace{-0.5cm}
\includegraphics[height=3.5cm,width=4.2cm]{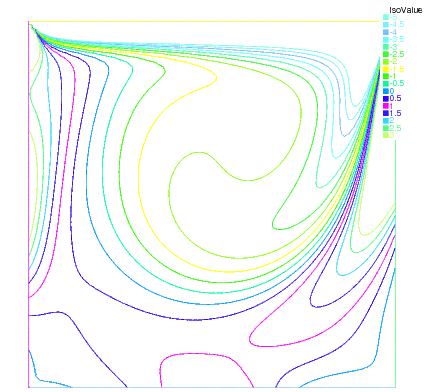}
\hspace{-0.5cm}
\includegraphics[height=3.5cm,width=4.2cm]{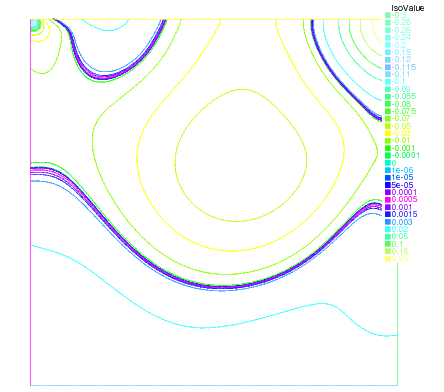}
\hspace{-0.5cm}
\includegraphics[height=3.5cm,width=4.2cm]{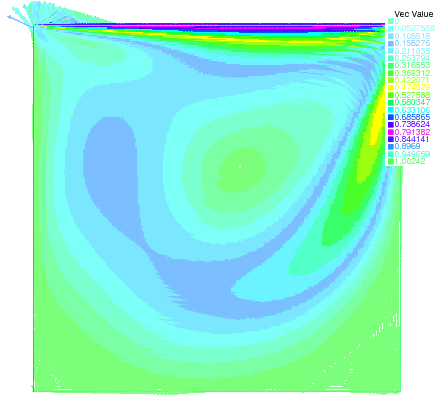}\\
\end{center}
\vspace{-0.3cm}
\caption{From top to bottom: Algorithms \ref{ALGO1G1},
\ref{NSALGO1}, \ref{NSALGO2}, \ref{NSALGO3} and \ref{NSALGO4}, and
from left to right~: flow, vorticity, pressure and velocity for
$\Delta t= 10^{-2}, Re=400, \tau=0.5, T=35, dim (X_H)=6561,
dim(Y_H)=1681, dim(X_h)=25921, dim(Y_h)=6561$.} \label{FigNSRe400}
\end{figure}
\end{center}
\newpage
\begin{center}
\begin{figure}[!h]
\begin{center}
\hspace{-0.4cm}
\includegraphics[height=3.5cm,width=4.2cm]{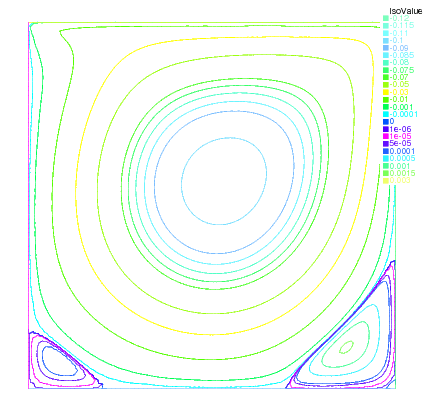}
\hspace{-0.5cm}
\includegraphics[height=3.5cm,width=4.2cm]{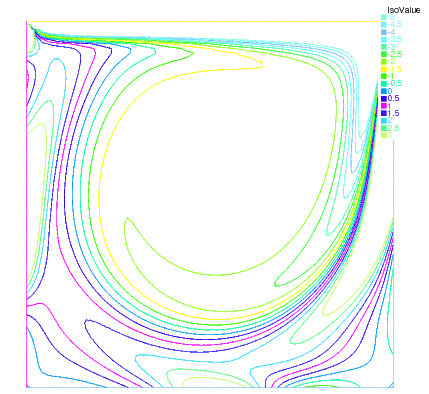}
\hspace{-0.5cm}
\includegraphics[height=3.5cm,width=4.2cm]{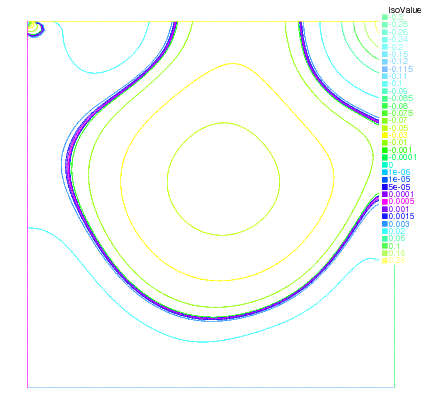}
\hspace{-0.5cm}
\includegraphics[height=3.5cm,width=4.2cm]{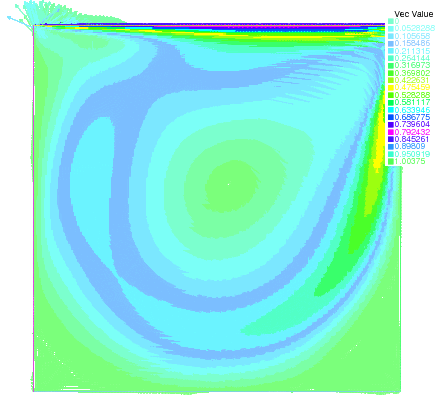}\\
%
\hspace{-0.4cm}
\includegraphics[height=3.5cm,width=4.2cm]{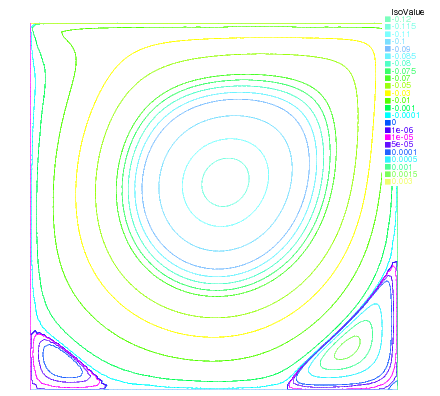}
\hspace{-0.5cm}
\includegraphics[height=3.5cm,width=4.2cm]{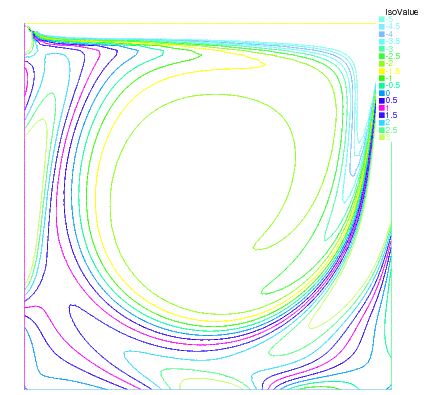}
\hspace{-0.5cm}
\includegraphics[height=3.5cm,width=4.2cm]{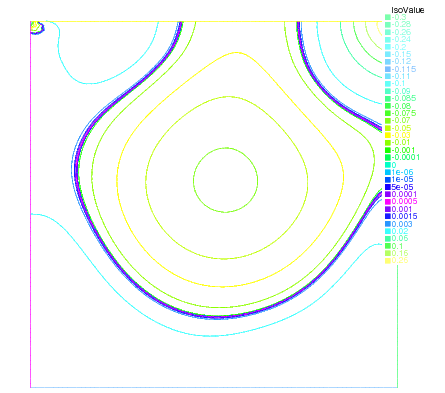}
\hspace{-0.5cm}
\includegraphics[height=3.5cm,width=4.2cm]{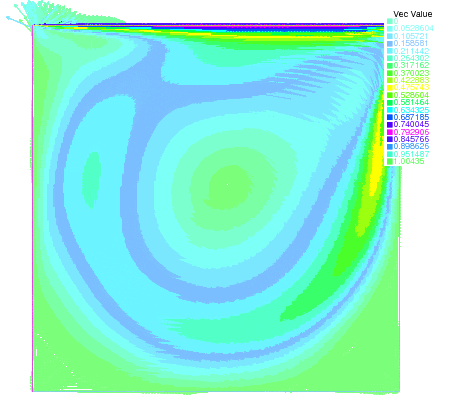}\\
%
\hspace{-0.4cm}
\includegraphics[height=3.5cm,width=4.2cm]{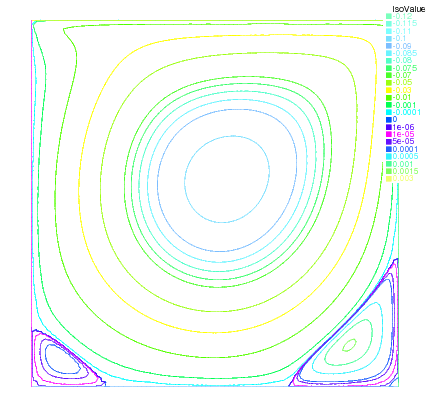}
\hspace{-0.5cm}
\includegraphics[height=3.5cm,width=4.2cm]{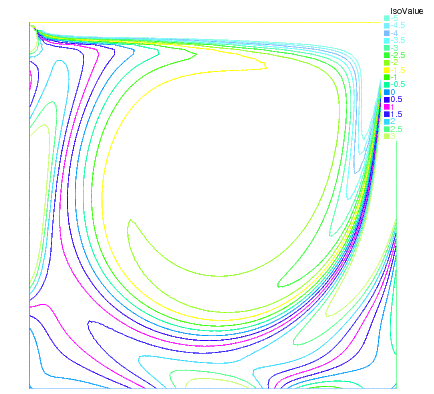}
\hspace{-0.5cm}
\includegraphics[height=3.5cm,width=4.2cm]{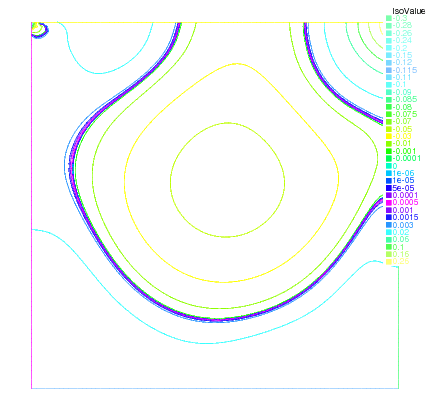}
\hspace{-0.5cm}
\includegraphics[height=3.5cm,width=4.2cm]{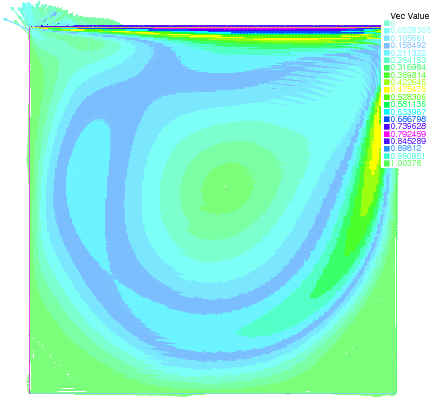}\\
%
\hspace{-0.4cm}
\includegraphics[height=3.5cm,width=4.2cm]{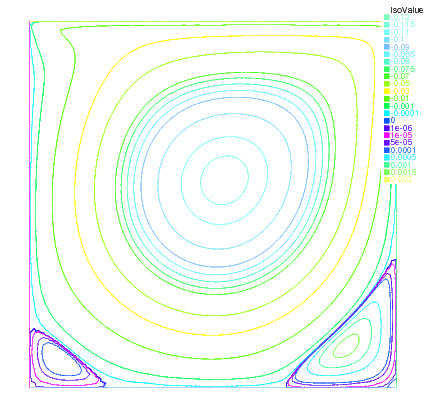}
\hspace{-0.5cm}
\includegraphics[height=3.5cm,width=4.2cm]{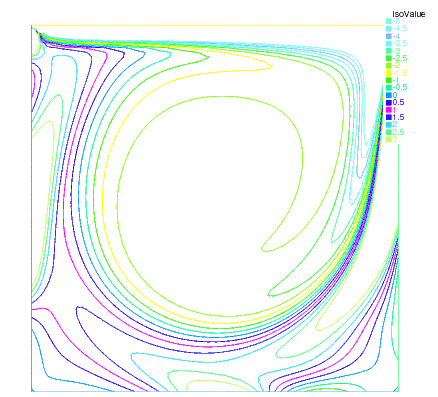}
\hspace{-0.5cm}
\includegraphics[height=3.5cm,width=4.2cm]{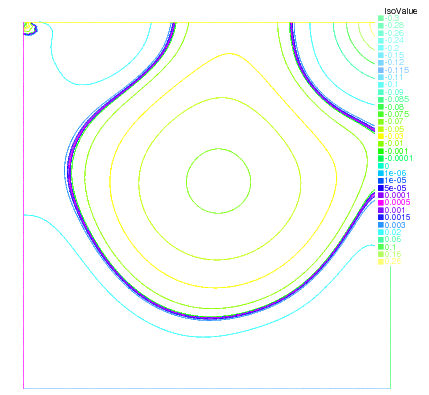}
\hspace{-0.5cm}
\includegraphics[height=3.5cm,width=4.2cm]{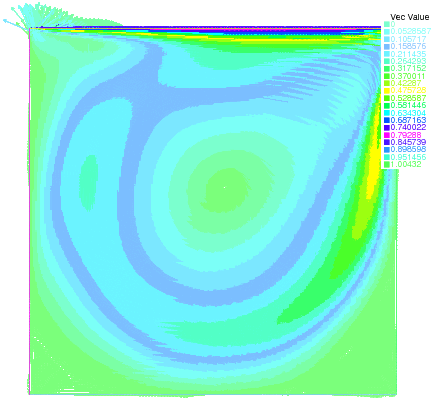}\\
\hspace{-0.4cm}
\includegraphics[height=3.5cm,width=4.2cm]{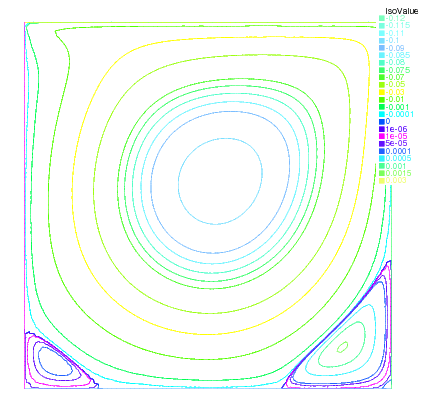}
\hspace{-0.5cm}
\includegraphics[height=3.5cm,width=4.2cm]{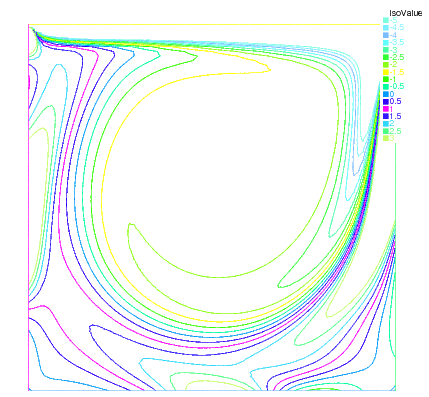}
\hspace{-0.5cm}
\includegraphics[height=3.5cm,width=4.2cm]{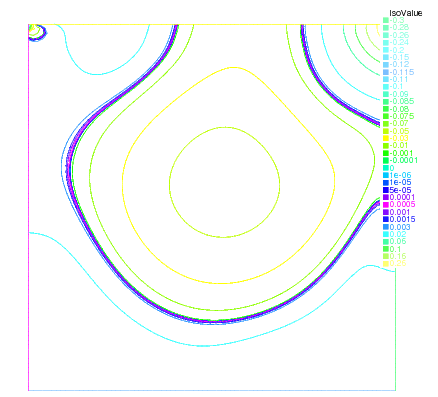}
\hspace{-0.5cm}
\includegraphics[height=3.5cm,width=4.2cm]{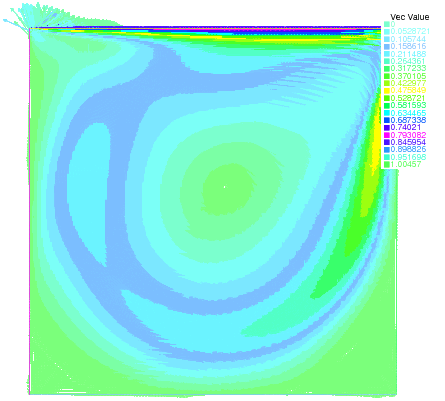}\\
\end{center}
\vspace{-0.3cm}
\caption{From top to bottom: Algorithms \ref{ALGO1G1},
\ref{NSALGO1}, \ref{NSALGO2}, \ref{NSALGO3} and \ref{NSALGO4}, and
from left to right~: flow, vorticity, pressure and velocity for
$\Delta t=5\times 10^{-3}, Re=1000, \tau=0.5, T=56, dim(X_H)=6561,
dim(Y_H)=1681, dim(X_h)=25921, dim(Y_h)=6561$.}
\label{FigNSRe1000}
\end{figure}
\end{center}
\newpage
\begin{center}
\begin{figure}[!h]
\begin{center}
\hspace{-0.4cm}
\includegraphics[height=3.5cm,width=4.2cm]{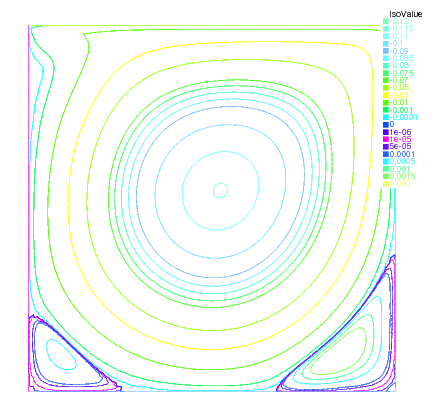}
\hspace{-0.5cm}
\includegraphics[height=3.5cm,width=4.2cm]{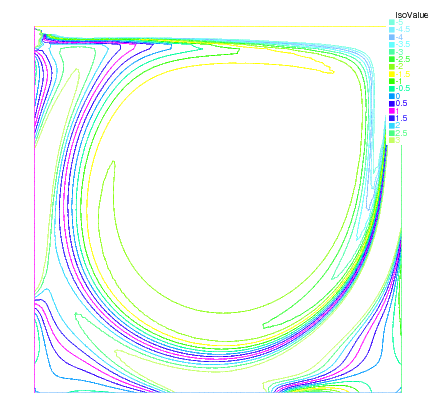}
\hspace{-0.5cm}
\includegraphics[height=3.5cm,width=4.2cm]{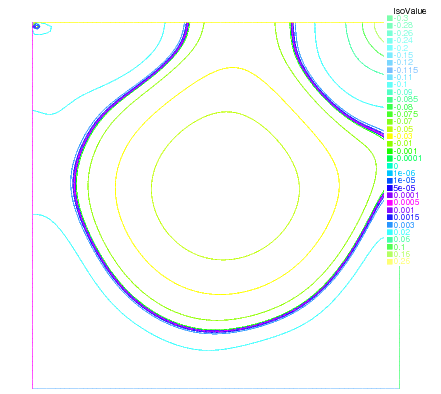}
\hspace{-0.5cm}
\includegraphics[height=3.5cm,width=4.2cm]{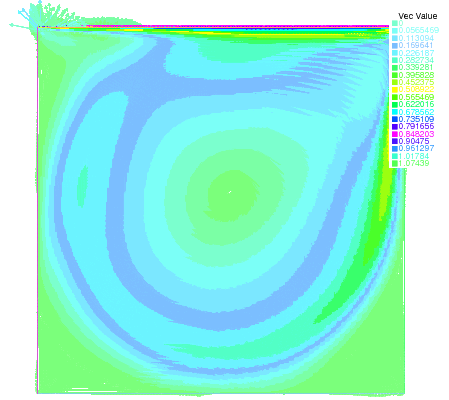}\\
%
\hspace{-0.4cm}
\includegraphics[height=3.5cm,width=4.2cm]{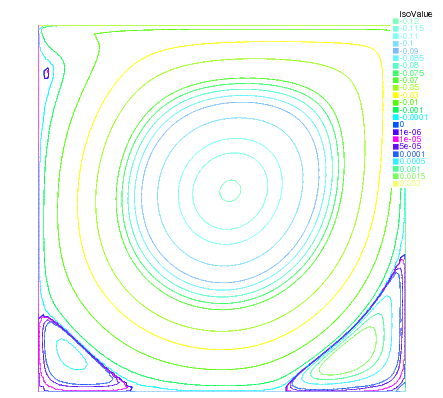}
\hspace{-0.5cm}
\includegraphics[height=3.5cm,width=4.2cm]{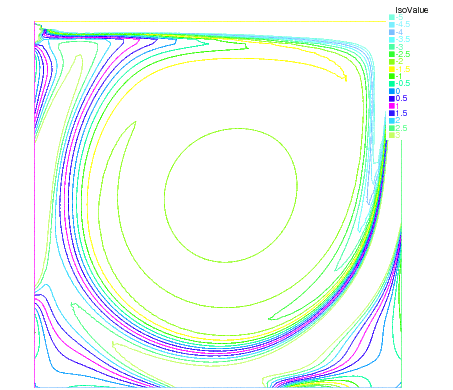}
\hspace{-0.5cm}
\includegraphics[height=3.5cm,width=4.2cm]{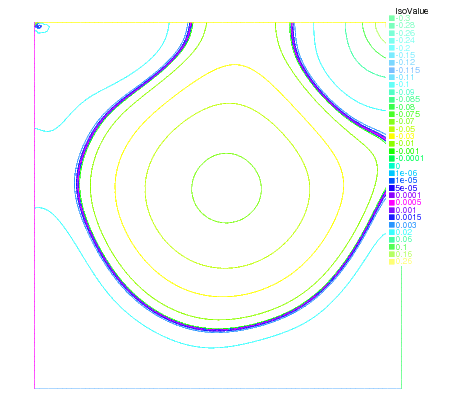}
\hspace{-0.5cm}
\includegraphics[height=3.5cm,width=4.2cm]{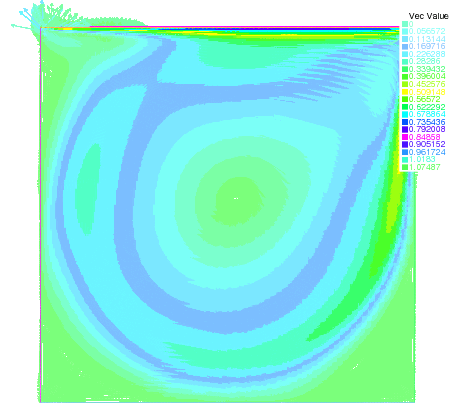}\\
%
\end{center}
%
\caption{From top to bottom: Algorithms \ref{ALGO1G1} and
\ref{NSALGO1} for $Re=2000, T=70$ respectively and from left to
right~: flow, vorticity, pressure and velocity for $\Delta t=
10^{-3}, \tau=0.5, dim(X_H)=6561, dim(Y_H)=1681, dim(X_h)=25921,
dim(Y_h)=6561$.} \label{FigNSRe2000}
\end{figure}
\end{center}
\begin{center}
\begin{figure}[!h]
\begin{center}
\hspace{-0.4cm}
\includegraphics[height=3.5cm,width=4.2cm]{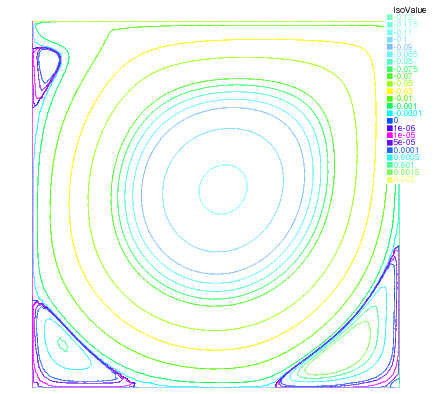}
\hspace{-0.5cm}
\includegraphics[height=3.5cm,width=4.2cm]{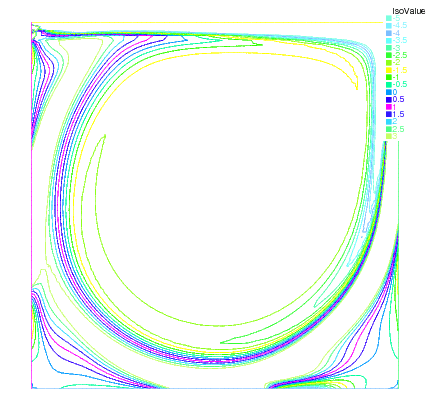}
\hspace{-0.5cm}
\includegraphics[height=3.5cm,width=4.2cm]{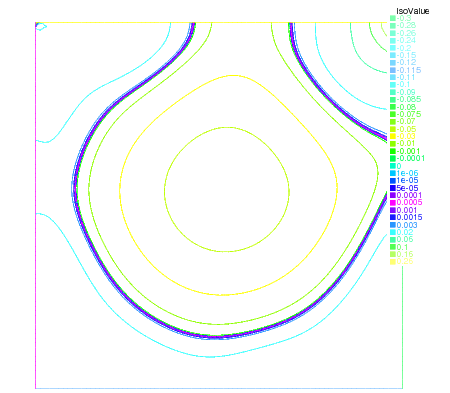}
\hspace{-0.5cm}
\includegraphics[height=3.5cm,width=4.2cm]{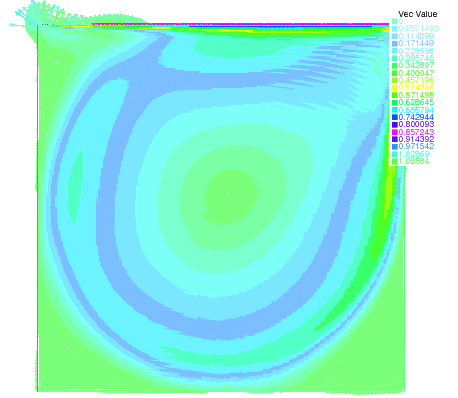}\\
\hspace{-0.4cm}
\includegraphics[height=3.5cm,width=4.2cm]{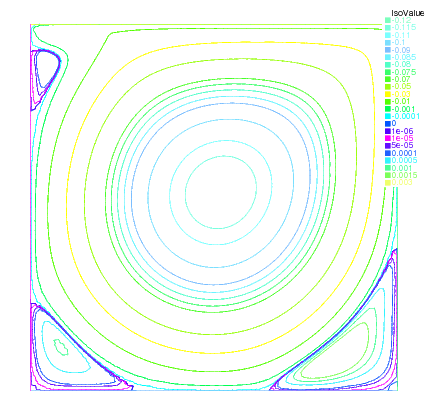}
\hspace{-0.5cm}
\includegraphics[height=3.5cm,width=4.2cm]{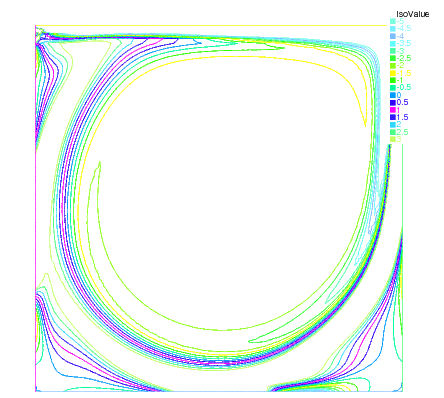}
\hspace{-0.5cm}
\includegraphics[height=3.5cm,width=4.2cm]{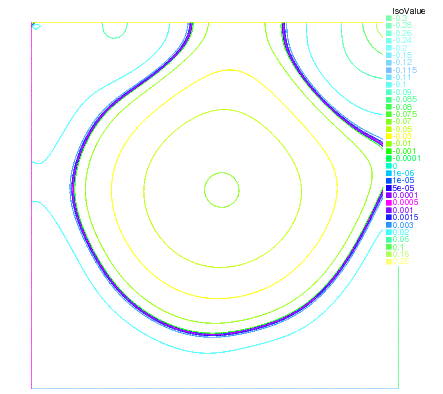}
\hspace{-0.5cm}
\includegraphics[height=3.5cm,width=4.2cm]{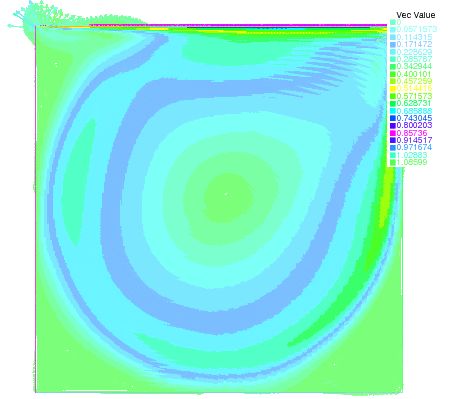}\\
\end{center}
%
\caption{From top to bottom: Algorithms \ref{ALGO1G1} and
\ref{NSALGO1} for $Re=3200, T=130$ respectively and from left to
right~: flow, vorticity, pressure and velocity for $\Delta t=
10^{-3}, \tau=0.5, dim(X_H)=6561, dim(Y_H)=1681, dim(X_h)=25921,
dim(Y_h)=6561$.} \label{FigNSRe3200}
\end{figure}
\end{center}
\newpage
\begin{center}
\begin{figure}[!h]
\begin{center}
\hspace{-0.4cm}
\includegraphics[height=3.5cm,width=4.2cm]{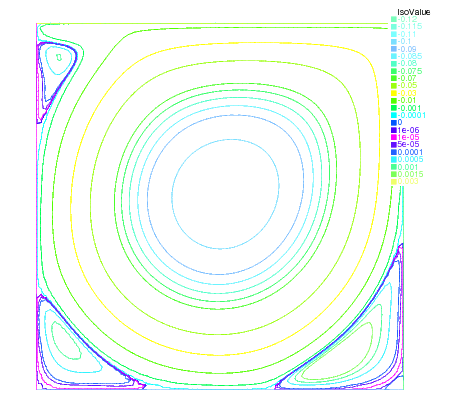}
\hspace{-0.5cm}
\includegraphics[height=3.5cm,width=4.2cm]{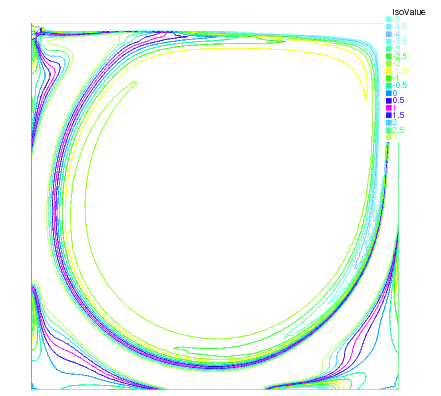}
\hspace{-0.5cm}
\includegraphics[height=3.5cm,width=4.2cm]{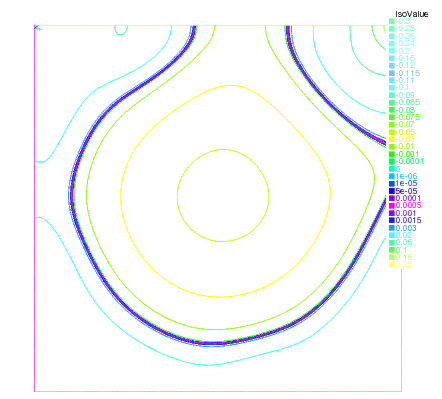}
\hspace{-0.5cm}
\includegraphics[height=3.5cm,width=4.2cm]{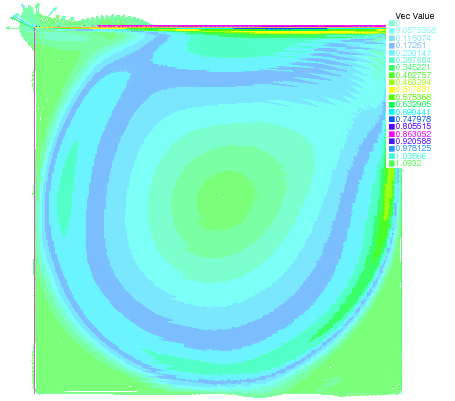}\\
\hspace{-0.4cm}
\includegraphics[height=3.5cm,width=4.2cm]{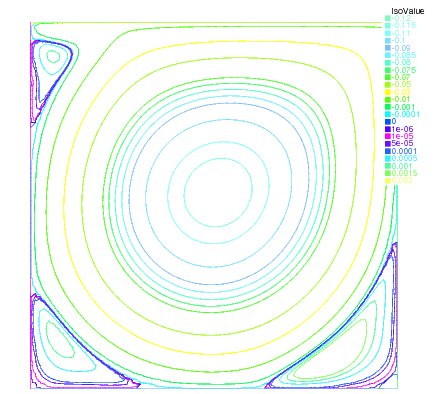}
\hspace{-0.5cm}
\includegraphics[height=3.5cm,width=4.2cm]{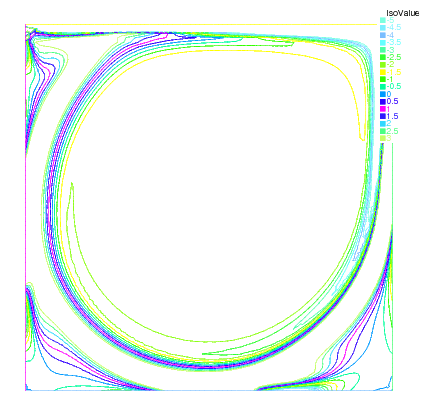}
\hspace{-0.5cm}
\includegraphics[height=3.5cm,width=4.2cm]{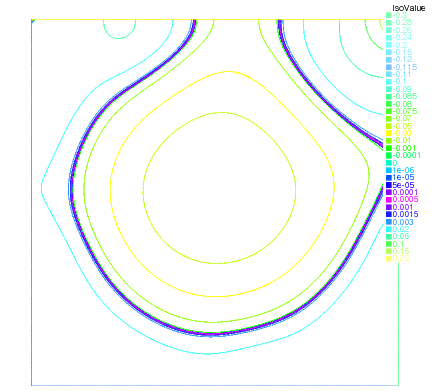}
\hspace{-0.5cm}
\includegraphics[height=3.5cm,width=4.2cm]{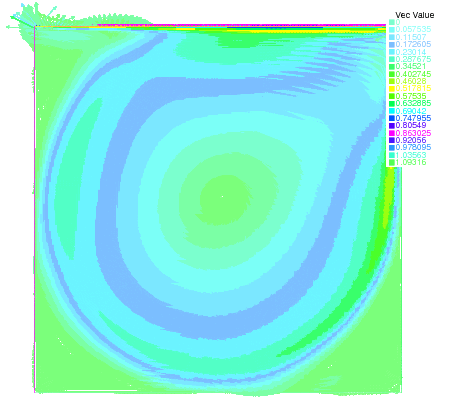}\\

\end{center}
%
\caption{ From top to bottom: Algorithm \ref{ALGO1G1} and
\ref{NSALGO1} for $Re=5000, T=148$ respectively and from left to
right: flow, vorticity, pressure and velocity for $\Delta t=
10^{-3},\tau=0.5, dim(X_H)=6561, dim(Y_H)=1681, dim(X_h)=25921,
dim(Y_h)=6561$.} \label{FigNSRe5000}
\end{figure}
\end{center}
\begin{center}
\begin{figure}[!h]
\begin{center}
\hspace{-0.4cm}
\includegraphics[height=3.5cm,width=4.2cm]{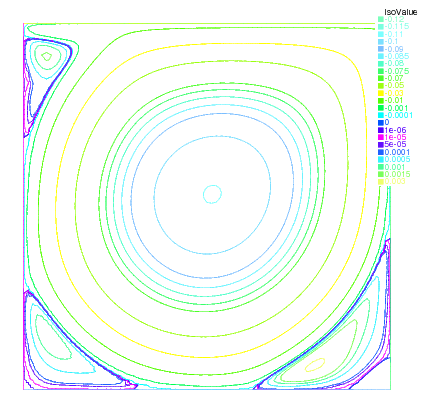}
\hspace{-0.5cm}
\includegraphics[height=3.5cm,width=4.2cm]{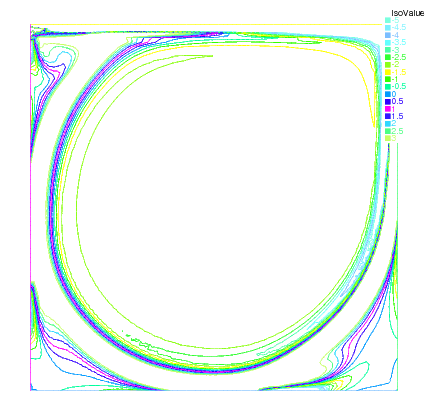}
\hspace{-0.5cm}
\includegraphics[height=3.5cm,width=4.2cm]{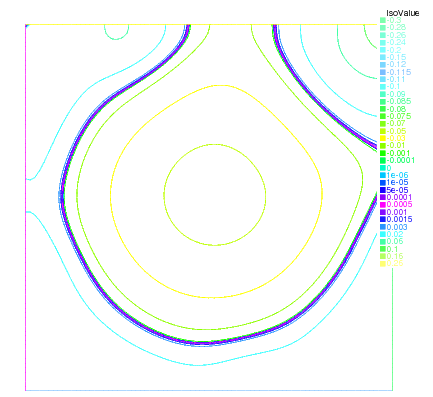}
\hspace{-0.5cm}
\includegraphics[height=3.5cm,width=4.2cm]{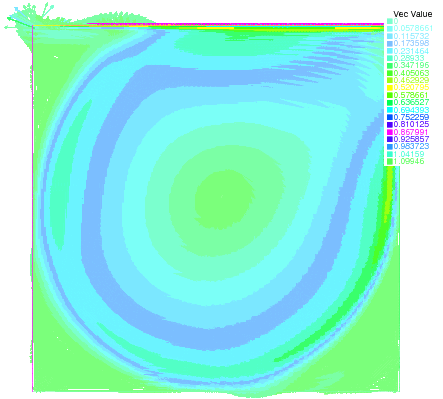}\\
\hspace{-0.4cm}
\includegraphics[height=3.5cm,width=4.2cm]{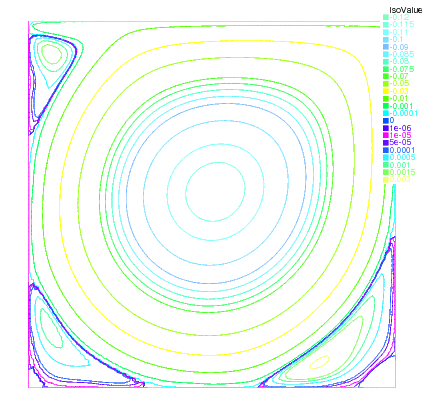}
\hspace{-0.5cm}
\includegraphics[height=3.5cm,width=4.2cm]{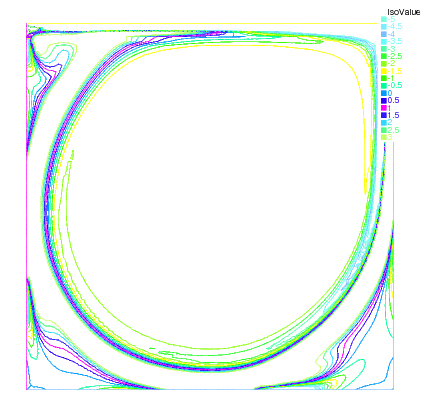}
\hspace{-0.5cm}
\includegraphics[height=3.5cm,width=4.2cm]{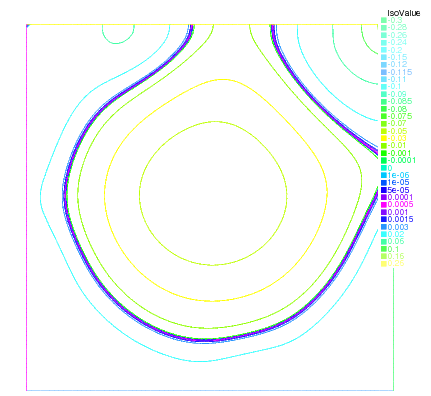}
\hspace{-0.5cm}
\includegraphics[height=3.5cm,width=4.2cm]{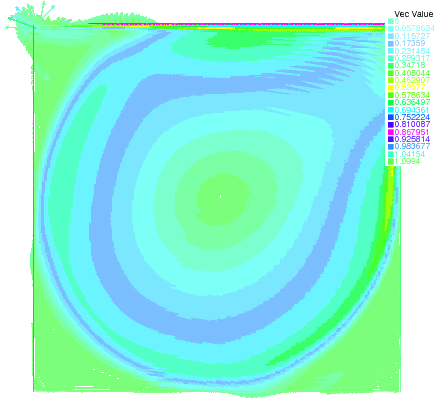}\\

\end{center}
%
\caption{ From top to bottom: Algorithms \ref{ALGO1G1} and
\ref{NSALGO1} for $Re=7500, T=195$ respectively and from left to
right~: flow, vorticity, pressure and velocity for $\Delta
t=5\times 10^{-4}, \tau=0.5, dim(X_H)=6561, dim(Y_H)=1681,
dim(X_h)=25921, dim(Y_h)=6561$.} \label{FigNSRe7500}
\end{figure}
\end{center}
\newpage

These flow patterns in Figures \ref{FigNSRe400},
\ref{FigNSRe1000}, \ref{FigNSRe2000}, \ref{FigNSRe3200},
\ref{FigNSRe5000}, \ref{FigNSRe7500} fit with
earlier results of Bruneau et al. \cite{Bruneau1990}, Ghia et al. \cite{Ghia1982}, Goyon \cite{Goyon1986}, Pascal \cite{Pascal1992}
and Vanka \cite{Vanka1986}. Also, the numerical values and the localization of the extrema of the vorticity and of the stream function are good agreement with the ones given in these references. 
\clearpage
\subsubsection{CPU Time reduction}
%
In Figures \ref{FigCPU4Algo} - \ref{FigCPUAlgo13Re} 
we observe that the bi-grid schemes are faster in computation time than reference the implicit scheme
(Algorithm \ref{ALGO1G1}) 
applied on the fine space $V_h$. However, the gain in CPU time is mainly obtained in the transient phase. Indeed, since a stationary solution is
computed, the reference scheme will
need only one nonlinear iteration at each the time step in the neighborhood of the steady state: it means that
only a linear system has then to been solved at each iteration, exactly as for the semi-implicit scheme. Therefore, in a neighborhood of the steady state, an iteration of the bi-grid method needs an additional
implicit iteration on the coarse grid and becomes then more expensive in CPU time.\\
\begin{center}
\begin{figure}[!h]
\vspace{-4.3cm}
\begin{center}
\hspace{-0.9cm}
\includegraphics[height=11cm,width=7cm]{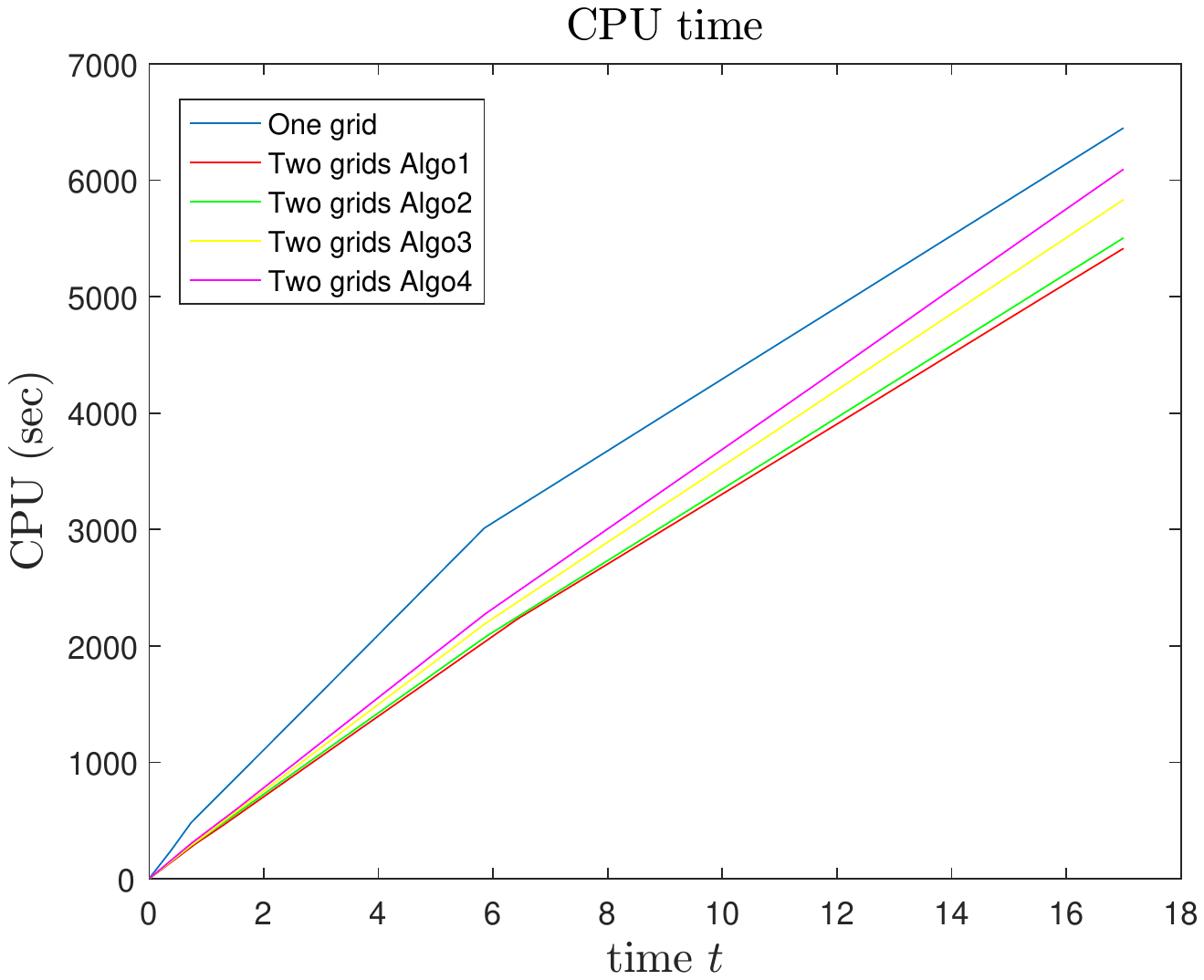}
\hspace{-2.5cm}
\includegraphics[height=11cm,width=7cm]{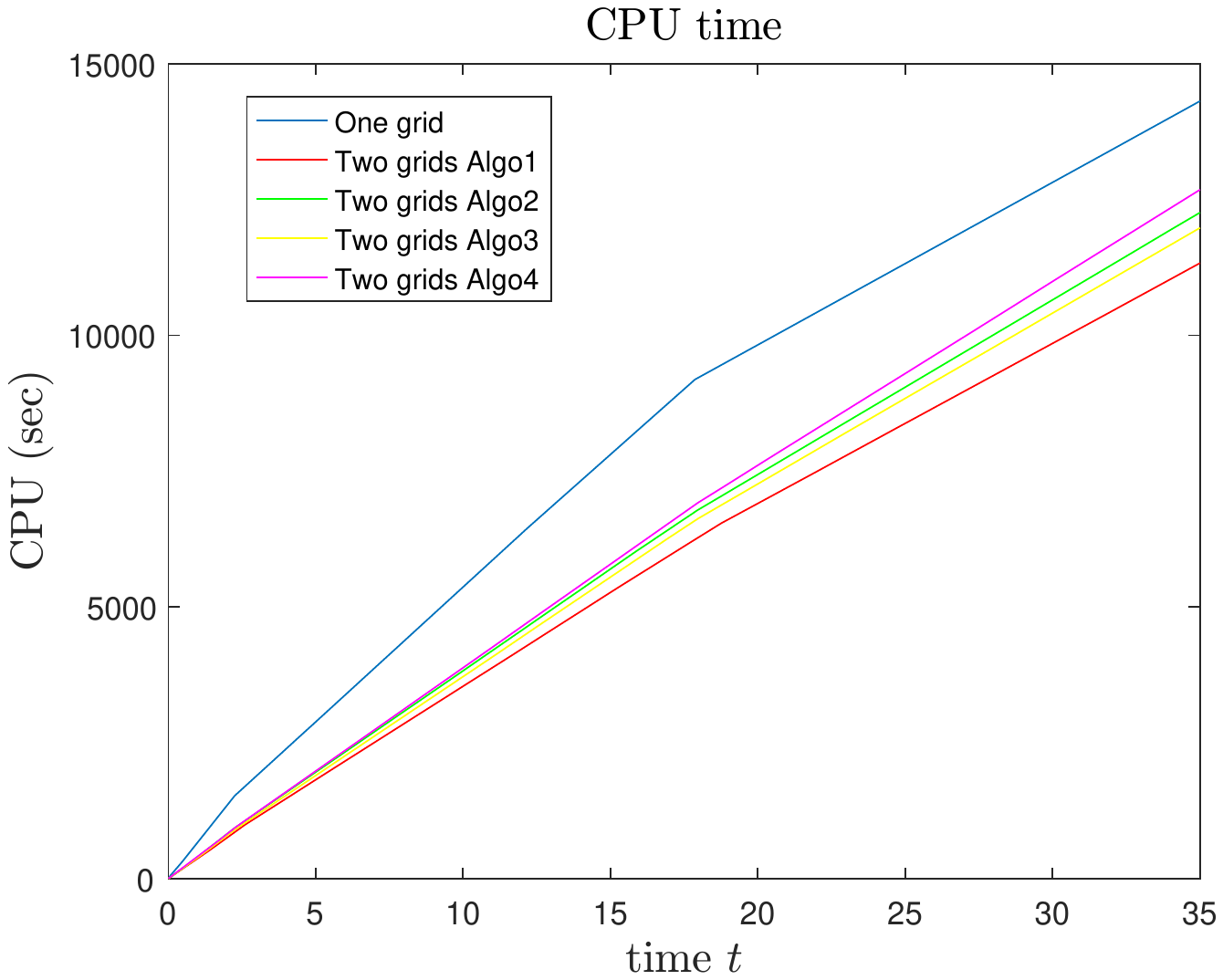}
\hspace{-2.5cm}
\includegraphics[height=11cm,width=7cm]{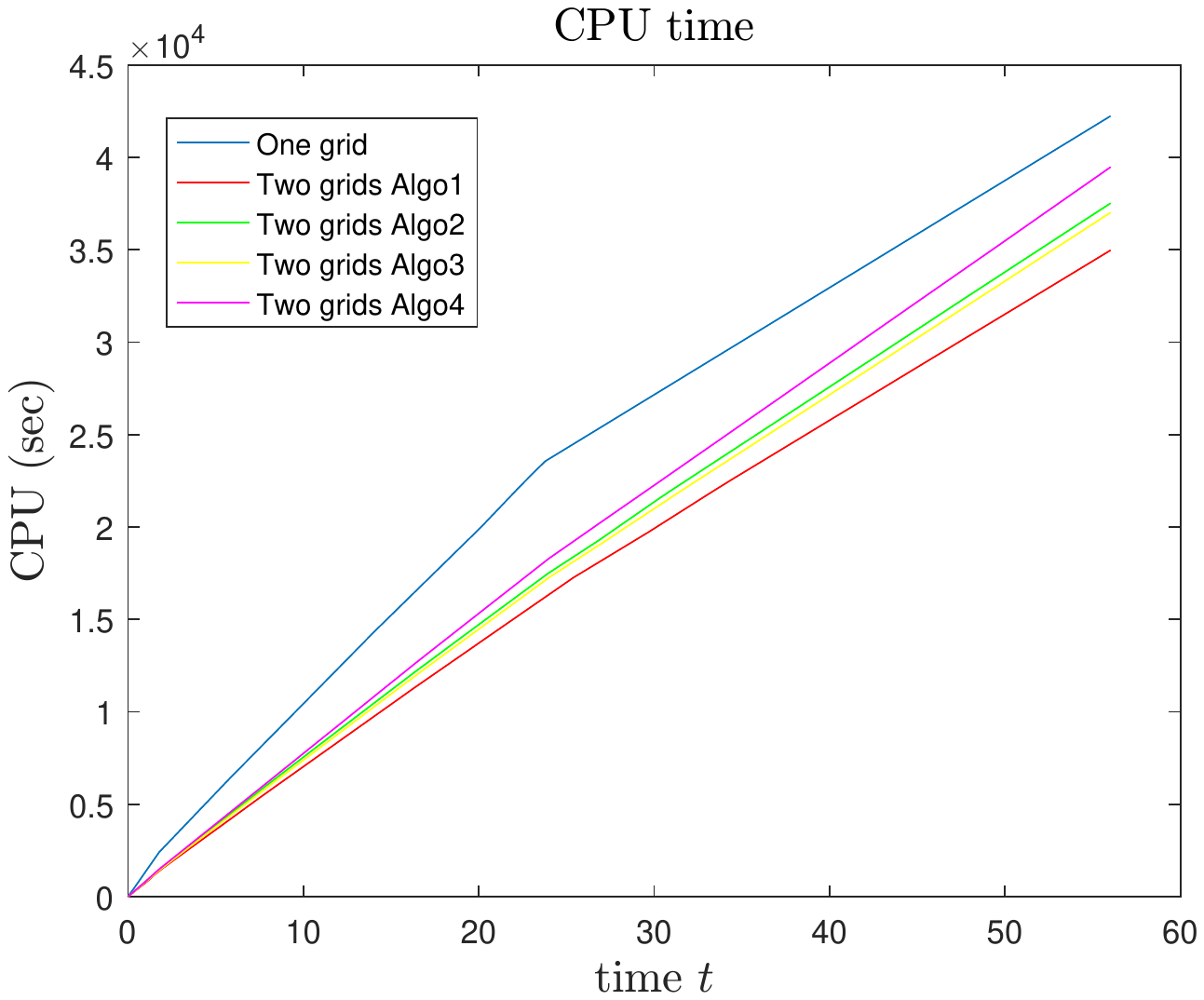}
\end{center}
\vspace{-3.7cm}
\caption{From left to right: CPU time of the scheme on a fine grid (Algorithm \ref{ALGO1G1}) and that of the scales
separation methods for $Re=100, 400, 1000$ respectively.} \label{FigCPU4Algo}
\end{figure}
\end{center}
\begin{center}
\begin{figure}[!h]
\vspace{-4.2cm}
\begin{center}
\hspace{-0.9cm}
\includegraphics[height=11cm,width=7cm]{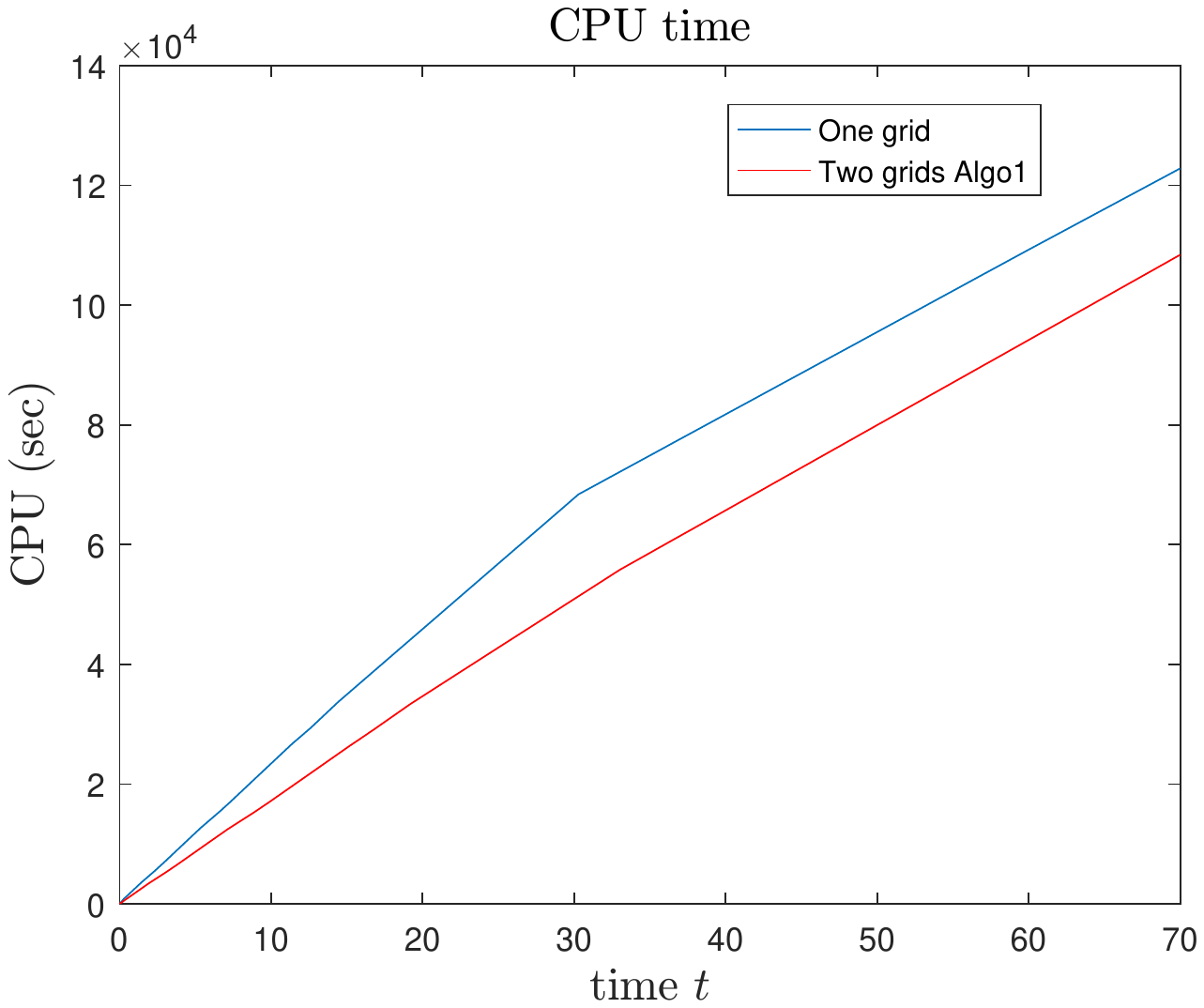}
\hspace{-2.5cm}
\includegraphics[height=11cm,width=7cm]{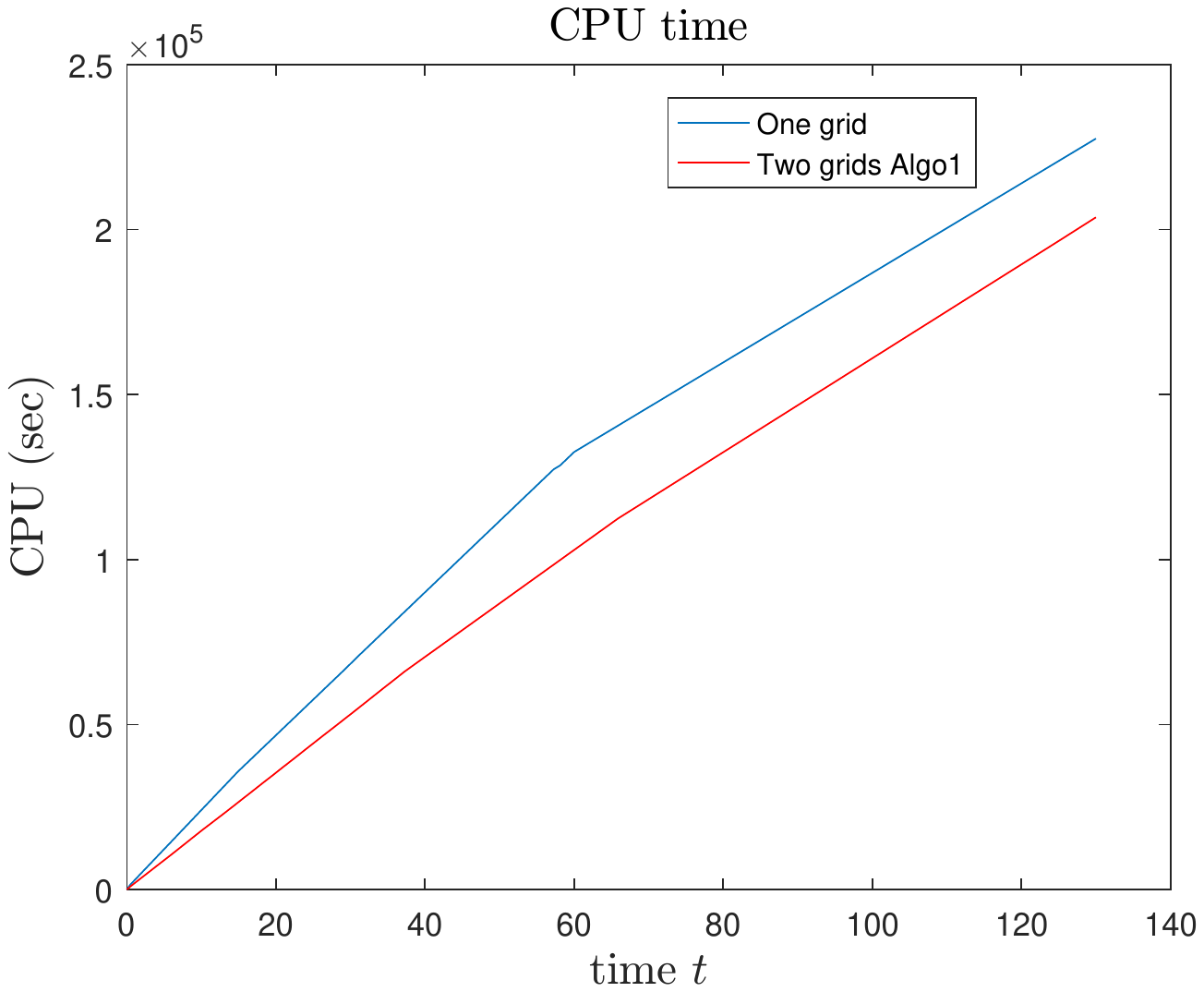}
\end{center}
\vspace{-3.4cm}
\caption{From left to right: CPU time of
the scheme on a fine grid (Algorithm \ref{ALGO1G1}) and that of the scales
separation method (Algorithm \ref{NSALGO1}) for $Re=2000, 3200$ respectively.}
\label{FigCPUAlgo13Re}
\end{figure}
\end{center}
\begin{center}
\begin{figure}[!h]
\vspace{-4.5cm}
\begin{center}
\hspace{-0.9cm}
\includegraphics[height=11cm,width=7cm]{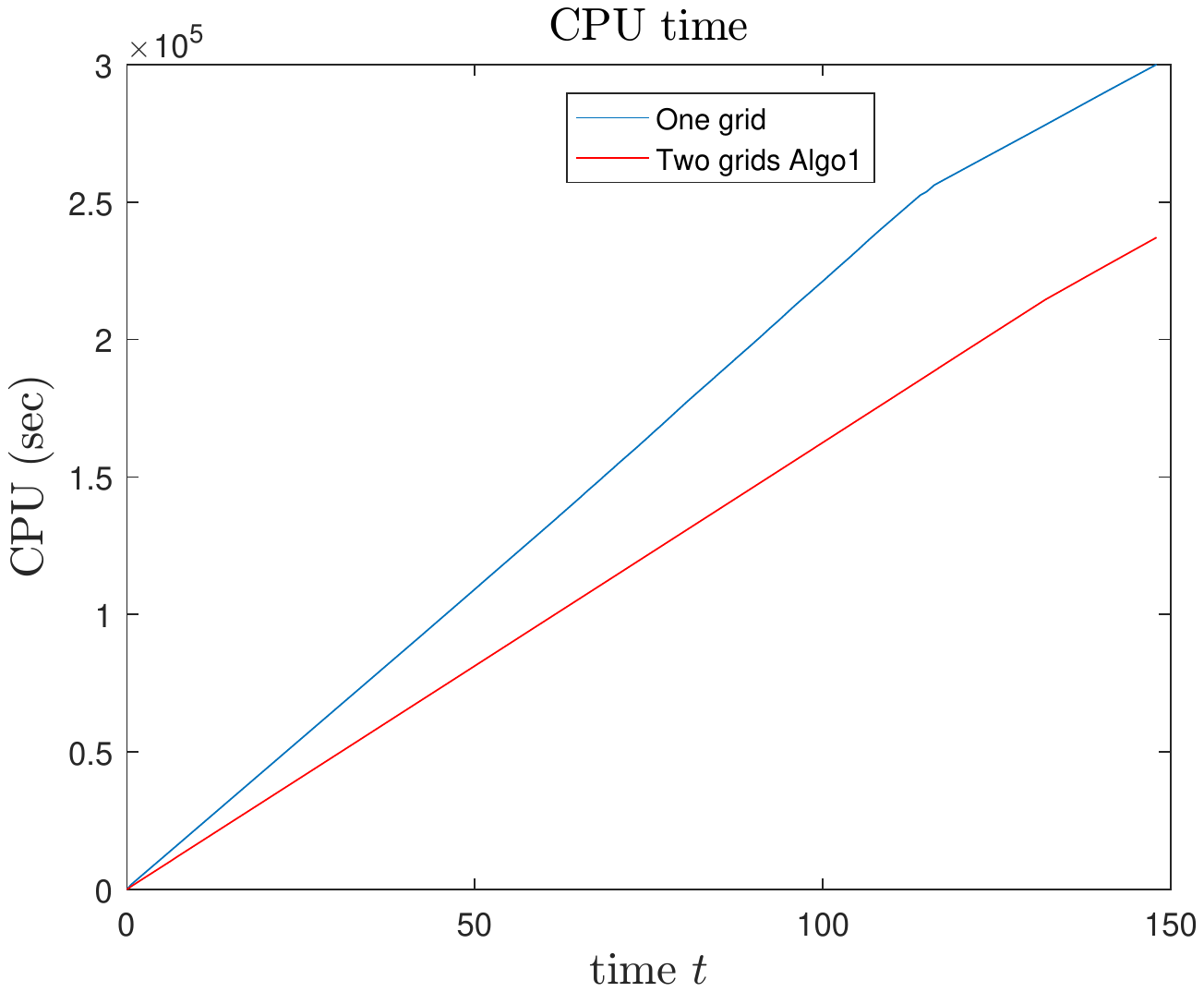}
\hspace{-2.5cm}
\includegraphics[height=11cm,width=7cm]{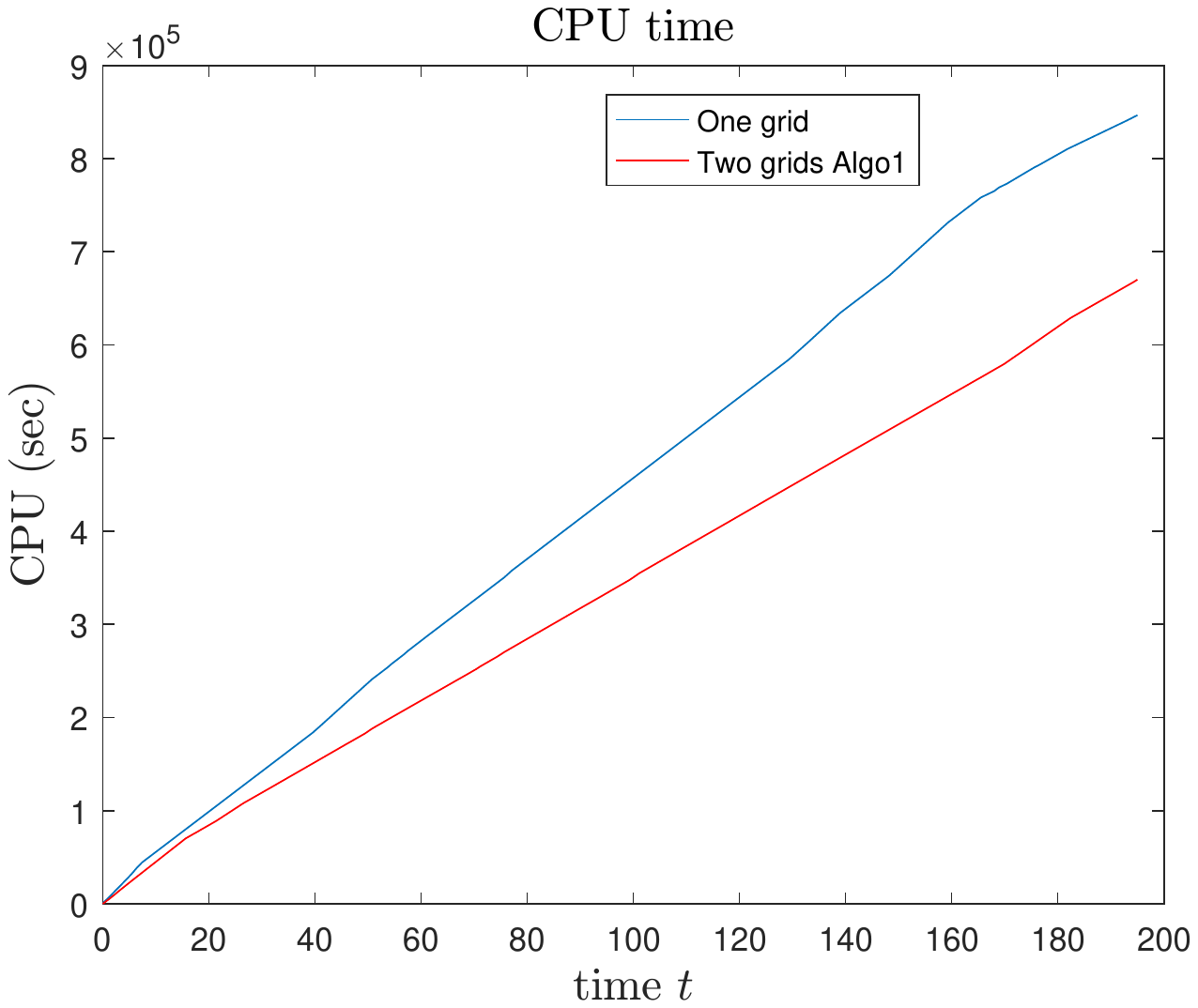}
\end{center}
\vspace{-3.4cm}
\caption{From left to right: CPU time of the
scheme on a fine grid (Algorithm \ref{ALGO1G1}) and that of the
scales separation method (Algorithm \ref{NSALGO1}) for $Re=5000, 7500$ respectively.}
\label{FigCPUAlgo13Re}
\end{figure}
\end{center}

We also compute the number of nonlinear iterations as
a function of time and compute the $L^2$ norm of
$\ds\frac{\partial u}{\partial t}$; we use the approximation $\ds\frac{\partial u}{\partial t}|_{t=t_k}\simeq \ds\frac{u^{k+1}-u^k}{\Delta t}$.%
We identify the time from which the reference scheme (Algorithm \ref{ALGO1G1}) reduces to only one linear iteration at each time step
and we define $\theta_s$ as the associated threshold value of
$\|\ds\frac{\partial u}{\partial t}\|_{L^2(\Omega}$. We then can obtain an heuristic criteria to define a simple strategy to save more CPU time and we modify the bi-grid scheme as follows:
as long as $\|\ds\frac{\partial u}{\partial t}\|>\theta_s$
we apply the bi-grid scheme and as soon as
$\|\ds\frac{\partial u}{\partial t}\|_{L^2(\Omega}\le \theta_s$ we apply the reference scheme on the fine space $V_h$. We
represent in Tables \ref{CPUtableRe100}, \ref{CPUtableRe400},
\ref{CPUtableRe1000}, \ref{CPUtablealgo1Re2000320050007500} the
calculation times for the four algorithms. We notice that the gain
ratio is now equal to $0.7$ instead of $0.9$ when no threshold strategy was applied on the values of $\ds\frac{\partial u}{\partial t}$ .
\begin{table}[!h]
\begin{center}
\begin{tabular}{|l||l||l||l||l|}
\hline
$80\times 40$ points &   Algorithm \ref{NSALGO1} &  Algorithm \ref{NSALGO2} &  Algorithm \ref{NSALGO3} &  Algorithm \ref{NSALGO4} \\
\hline 
\hline
$t_{1G}$ at $T=17$(in sec.) & 6448.39 & 6448.39 & 6448.39 &  6448.39 \\
\hline
\hline
CPU (in sec.) & 5414 & 5504.2 & 5833.06 & 6094.01 \\
\hline
\hline
$\ds\frac{t_{2G} }{t_{1G}}$  & 0.8396 & 0.8536 & 0.9046 & 0.9450 \\
\hline
 \hline
$\theta_s$ & 0.000528406 & 0.000528406 & 0.000528406 & 0.000528406 \\
%
$t_{1G}$ (in sec.) & 3015.92 & 3015.92 & 3015.92 & 3015.92 \\
\hline

\hline
CPU (in sec.) & 2036.78 & 2074.37 & 2190.44 & 2273.63  \\
\hline
\hline
$\ds\frac{t_{2G} }{t_{1G}}$  & 0.675 & 0.687 & 0.726 & 0.784\\
\hline
\end{tabular}
\end{center}
\vspace{-0.5cm}
\caption{Non incremental Navier-Stokes
$Re=100,\Delta t=10^{-2}$ and $\mathbb{P}_2/\mathbb{P}_1$.}\label{CPUtableRe100}
\end{table}
\begin{table}[!h]
\begin{center}
\begin{tabular}{|l||l||l||l||l|}
\hline
$80\times 40$ points &  Algorithm \ref{NSALGO1} &  Algorithm \ref{NSALGO2} &  Algorithm \ref{NSALGO3} &  Algorithm \ref{NSALGO4} \\
%
\hline 
\hline
$t_{1G}$ at $T=35$(in sec.) & 14311.1 & 14311.1 & 14311.1 &  14311.1 \\
\hline
\hline
CPU (in sec.) & 11330.3 & 12259.1 & 11978 & 12680.4 \\
\hline
\hline
$\ds\frac{t_{2G} }{t_{1G}}$  & 0.7917 & 0.8566 & 0.8369 & 0.8860 \\
\hline
 \hline
$\theta_s$ & 0.000438706 & 0.000438706 & 0.000438706 & 0.000438706\\
%
$t_{1G}$ (in sec.) & 9191.99  &  9191.99 & 9191.99 &  9191.99\\
\hline

\hline
CPU (in sec.) & 6247.12 & 6756.45 & 6595.52 & 6886.03 \\
\hline
\hline
$\ds\frac{t_{2G} }{t_{1G}}$  & 0.679 & 0.735 & 0.717 & 0.749 \\
%
\hline
\end{tabular}
\end{center}
\caption{Non incremental Navier-Stokes
$Re=400,\Delta t=10^{-2}$ and $\mathbb{P}_2/\mathbb{P}_1$.}\label{CPUtableRe400}
\end{table}
\begin{table}[!h]
\begin{center}
\begin{tabular}{|l||l||l||l||l|}
\hline
$80\times 40$ points &  Algorithm \ref{NSALGO1} &  Algorithm \ref{NSALGO2} &  Algorithm \ref{NSALGO3} &  Algorithm \ref{NSALGO4} \\
\hline 
\hline
$t_{1G}$ at $T=56$(in sec.) & 42230.3 & 42230.3 & 42230.3 &  42230.3\\
\hline
\hline
CPU (in sec.) & 34967.7 & 37511.5 & 37009.8 & 39466.6 \\
\hline
\hline
$\ds\frac{t_{2G} }{t_{1G}}$  & 0.8280 & 0.8882 & 0.8764 & 0.9346 \\
%
\hline
 \hline
$\theta_s$ & 0.000948096 & 0.000948096 & 0.000948096 & 0.000948096\\
%
$t_{1G}$ (in sec.) & 23564.3 & 23564.3 & 23564.3 & 23564.3 \\
\hline

\hline
CPU (in sec.) & 16212.7 & 17382 & 17121.7 & 18146.9 \\
\hline
\hline
$\ds\frac{t_{2G} }{t_{1G}}$  & 0.688 & 0.737 & 0.726 & 0.77 \\
\hline
\end{tabular}
\end{center}
\caption{Non incremental Navier-Stokes $Re=1000,\Delta t=5\times 10^{-3}$ and $\mathbb{P}_2/\mathbb{P}_1$.}\label{CPUtableRe1000}
\end{table}
\begin{table}[!h]
\begin{center}
\begin{tabular}{|l||l||l||l||l|}
\hline
$80\times 40$ points &  $Re=2000$     &  $Re=3200$ &  $Re=5000$ & $Re=7500$ \\
       & $\Delta t=10^{-3}$  & $\Delta t=10^{-3}$  & $\Delta t=10^{-3}$   &   $\Delta t=5\times10^{-4}$\\
       \hline
 Final time & $T=70$ & $T=130$  & $T=148$ & $T=195$\\
\hline 
\hline
$t_{1G}$ at $T$(in sec.) & 122826 & 227543 & 299998 & 846647 \\
\hline
\hline
CPU (in sec.) &108401 & 203692 & 237105 &  670244\\
\hline
\hline
$\ds\frac{t_{2G} }{t_{1G}}$  & 0.8826 & 0.8952 & 0.7904 & 0.791 \\
\hline
 \hline
$\theta_s$ & 0.00521845 & 0.00502797 & 0.00321512 & 0.00462578\\
%
$t_{1G}$ (in sec.) & 68396.7 & 132604 & 256175 &  758398\\
\hline

\hline
CPU (in sec.) & 51375.2 & 102922 &  188600&  565283\\
\hline
\hline
$\ds\frac{t_{2G} }{t_{1G}}$  & 0.751 & 0.776 & 0.736 & 0.745  \\
\hline
\end{tabular}
\end{center}
\vspace{-0.5cm} \caption{Non incremental Navier-Stokes for $Re=2000, 3200, 5000, 7500$ and $\mathbb{P}_2/\mathbb{P}_1$
}\label{CPUtablealgo1Re2000320050007500}
\end{table}
\pagebreak
\subsubsection{Enhanced Stability : comparison with a semi-implicit scheme}
As presented in the previous sections, the bi-grid methods are faster than the unconditionally stable reference scheme (Algorithm \ref{ALGO1G1}).
To illustrate the stabilization properties of the bi-grid schemes, we make now a comparison with the semi-implicit scheme (applied on whole the fine space $V_h$) and which consists in treating the nonlinear term $(u^k\cdot\nabla) u^k$ explicitly.
Particularly we compare the maximum time steps $\Delta t$ that can be taken for each of the schemes and compare the respective final  CPU time on the time interval $[0,T_f]$.
%
\noindent We  give hereafter in Tables \ref{StabRe100}, \ref{StabRe400} and \ref{StabRe1000} the results comparing the stability of the bi-grid Algorithms \ref{NSALGO1}, \ref{NSALGO2}, \ref{NSALGO3}, \ref{NSALGO4}  and the following semi-implicit
scheme:
\begin{center}
\begin{minipage}[H]{16cm}
  \begin{algorithm}[H]
    \caption{Non-incremental semi-implicit scheme}\label{ALGO1G3}
    \begin{algorithmic}[1]
     \For{$k=0,1, \cdots$}
           \State {\bf Find $u_h^\ast$ in $X_h$}
           ($\ds\frac{u_h^\ast-u_h^k}{\Delta t},\psi_h) + \nu(\nabla
u_h^\ast,\nabla \psi_h)+((u_h^k\cdot\nabla)u_h^k,\psi_h)
=(f,\psi_h), \forall \psi_h \in X_h$
\State {\bf Find $p_h^{k+1}$ in $Y_h$} $(div\;\nabla p_h^{k+1},\chi_h)= (\ds\frac{div\;u_h^\ast}{\Delta
t},\chi_h), \forall \chi_h \in Y_h$
\State {\bf Find $u_h^{k+1}$ in $X_h$} ($u_h^{k+1}-u_h^\ast+\Delta t \nabla p_h^{k+1},\psi_h)=0, \forall
\psi_h \in X_h$
            \EndFor
 \end{algorithmic}
  \end{algorithm}
\end{minipage}
\end{center}

\begin{table}[!h]
\begin{center}
\begin{tabular}{|l|l|l|l|l|l|l|l|}
\hline
Scheme                 & $\tau$ & $\Delta t$ & Stability & $\Delta t$ & Stability & $\Delta t$ & Stability\\
\hline \hline
Algorithm \ref{NSALGO1} & $0.5$ & $0.01$ &  yes  & $0.05$ &  yes & $0.1$ &  yes \\
\hline
Algorithm \ref{NSALGO2} & $0.5$ & $0.01$ & yes  & $0.05$ &  yes & $0.1$ &  yes \\
\hline
Algorithm \ref{NSALGO3} & $0.5$ & $0.01$ &  yes  & $0.05$ &  yes  & $0.1$ &  yes \\
\hline
Algorithm \ref{NSALGO4} & $0.5$ & $0.01$ & yes  & $0.05$ &  yes  & $0.1$ &  yes \\
\hline
Algorithm \ref{ALGO1G3}      &     & $0.01$ & yes        &  $0.05$ &  yes     & $0.1$ &  yes \\
\hline
\end{tabular}
\end{center}
\vskip -0.5cm
\caption{$Re=100, dim(X_H)=6561, dim(Y_H)=1681, dim(X_h)=25921,
dim(Y_h)=6561, \mathbb{P}_2/\mathbb{P}_1$
elements.}\label{StabRe100}
\end{table}
%
%
%
\begin{table}[!h]
\begin{center}
\begin{tabular}{|l|l|l|l|l|l|l|l|l|l|}
\hline
Scheme                 & $\tau$ & $\Delta t$ & Stability & $\tau$& $\Delta t$ & Stability &$\tau$& $\Delta t$ & Stability\\
\hline \hline
Algorithm \ref{NSALGO1} & $0.5$ & $0.01$ &  yes & $0.5$ & $0.05$ &  no & $0.5$ & $0.1$ &  no \\
\cline{2-10}
& $-$       &$-$&$-$ & $30$ & $0.05$ &  yes  & $30$& $0.1$ &  yes \\
\hline \hline
Algorithm \ref{NSALGO2} & $0.5$ & $0.01$ & yes & $0.5 $ & $0.05$ &  yes & $0.5$ & $0.1$ &  yes \\
\cline{2-10}
 &$-$      &$-$&$-$& $30$ & $0.05$  &  yes & $30$ & $0.1$ &  yes \\
\hline \hline
Algorithm \ref{NSALGO3} & $0.5$ & $0.01$ &  yes & $0.5$ & $0.05$ &  no & $0.5$ & $0.1$ &  no \\
\cline{2-10}
&$-$ &$-$&$-$& $30$ & $0.05$ &  yes & $35$ &  $0.1$ &  yes \\
\hline\hline
Algorithm \ref{NSALGO4} & $0.5$ & $0.01$ & yes & $0.5$ & $0.05$ &  yes & $0.5$ & $0.1$ &  yes \\
\cline{2-10}
&$-$        &$-$&$-$& $30$  & $0.05$ &  yes & $30$ &  $0.1$ &  yes \\
\hline\hline
Algorithm \ref{ALGO1G3}           &       & $0.01$ & yes  & & $0.05$ &  no & &  $0.1$ &  no \\
\hline
\end{tabular}
\end{center}
\caption{$Re=400, dim(X_H)=6561, dim(Y_H)=1681, dim(X_h)=25921,
dim(Y_h)=6561, \mathbb{P}_2/\mathbb{P}_1$
elements.}\label{StabRe400}
\end{table}
%
%
%
\begin{table}[!h]
\begin{center}
\begin{tabular}{|l||l|l|l|l|l||l|l|l|}
\hline
Scheme                 & $\tau$ & $\Delta t$ & Stability & $\Delta t$ & Stability & $\tau$ & $\Delta t$ & Stability\\
\hline \hline
Algorithm \ref{NSALGO1} & $0.5$ & $0.005$ &  yes  & $0.007$ &  no & $100$ & $0.05$ &  yes\\
\hline
Algorithm \ref{NSALGO2} & $0.5$ & $0.005$ & yes  & $0.007$ &  yes & $0.5$ & $0.05$ &  yes \\
\hline
Algorithm \ref{NSALGO3} & $0.5$ & $0.005$ &  yes  & $0.007$ &  no & $150$ & $0.05$ &  yes \\
\hline
Algorithm \ref{NSALGO4} & $0.5$ & $0.005$ & yes   & $0.007$ &  yes & $ 0.5$ & $0.05$ &  yes\\
\hline
Algorithm \ref{ALGO1G3}      &       & $0.005$ & yes       &  $0.007$ & no &                & $0.05$ &  no \\
\hline
\end{tabular}
\end{center}
\caption{$Re=1000, dim(X_H)=6561, dim(Y_H)=1681, dim(X_h)=25921,
dim(Y_h)=6561, \mathbb{P}_2/\mathbb{P}_1$
elements.}\label{StabRe1000}
\end{table}
%
%
\noindent  We observe in Figure \ref{FigNSstabRe4001000} that the dynamics of the convergence to the steady state
of the bi-grid method is similar to the one of the implicit and
semi-implicit schemes. We note also that the stability of the bi-grid scheme
is guaranteed for large values of $\tau$. For example, for
$Re=1000$ and $\Delta t=0.007$, Algorithm \ref{NSALGO1} was
unstable for $\tau=0.5$ but taking a $\tau=30$ we can have a
stable scheme without deteriorating the  history of the convergence to the steady state.
In order to locate our gain at time step
$\Delta t$, we present in Table \ref{deltamax} a comparison
between the maximum time step allowing the stability of the
semi-implicit scheme and our Algorithm \ref{NSALGO1} for the
minimum value of $\tau$ necessary for stabilization.
\begin{center}
\begin{figure}[!h]
\begin{center}
\vspace{-3.5cm}
\includegraphics[height=13.7cm,width=8.5cm]{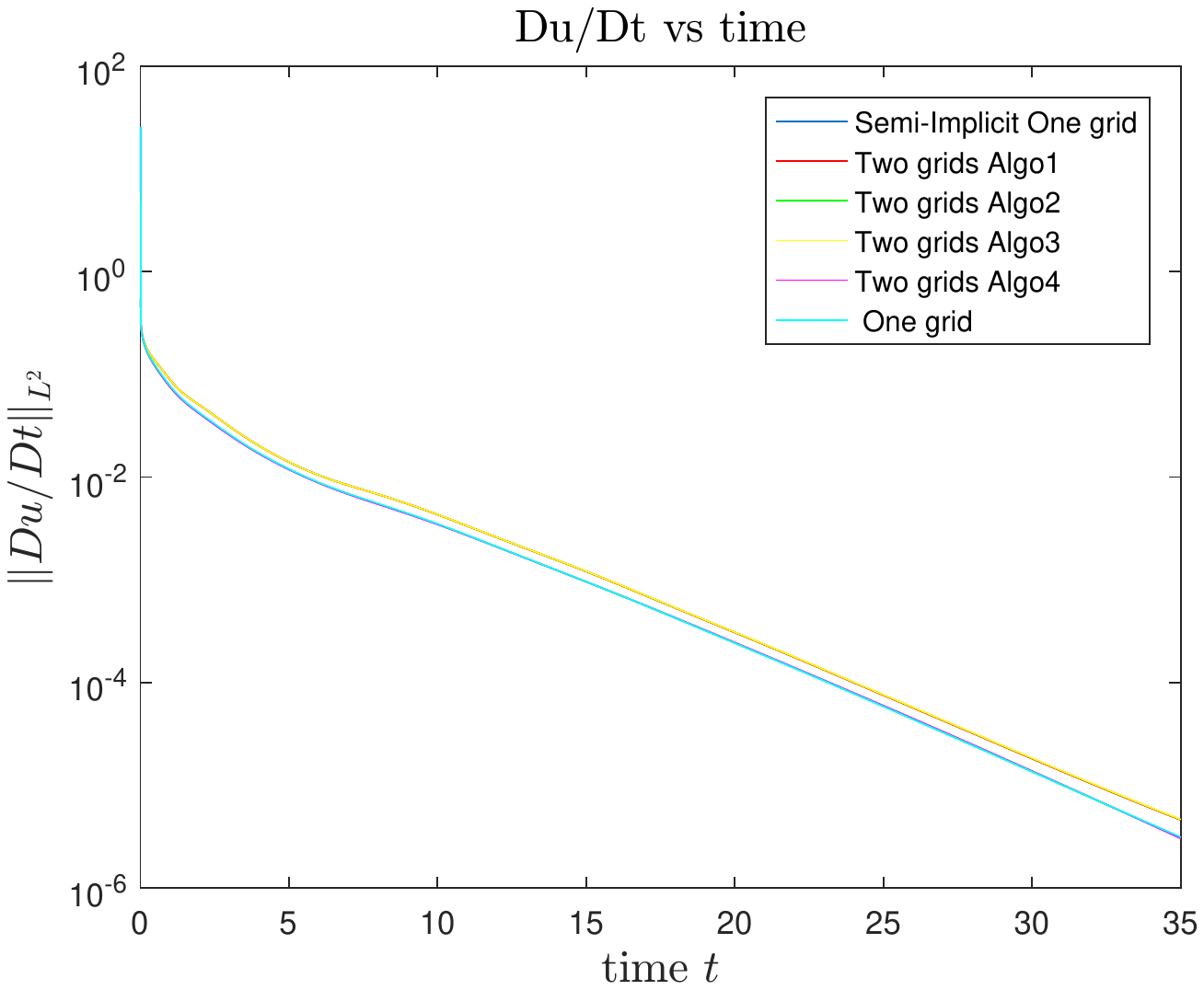}
\hspace{-1.7cm}
\includegraphics[height=13.7cm,width=8.5cm]{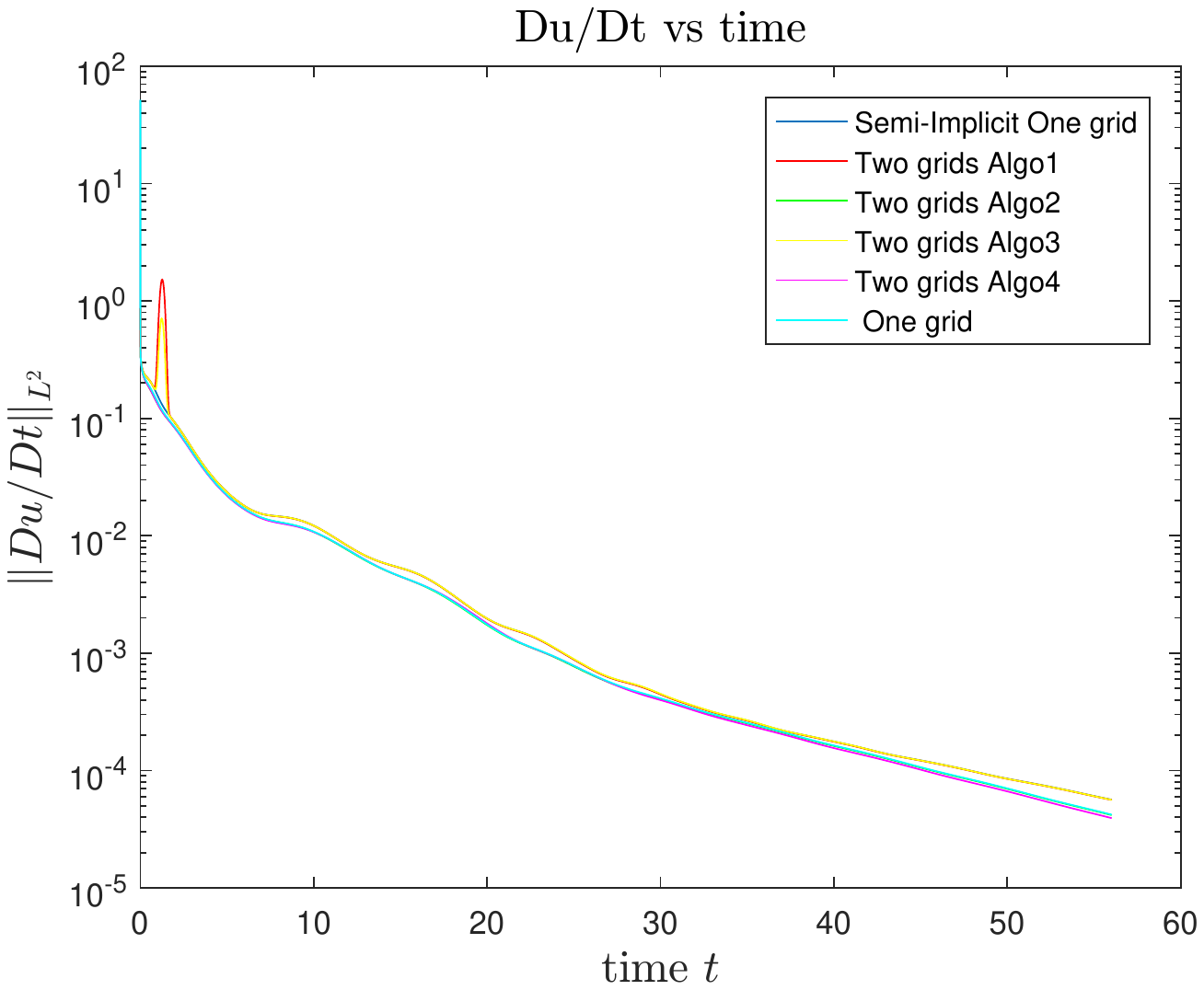}\\
\vspace{-3.9cm}
\end{center}
%
\caption{$\|\ds\frac{\partial u}{\partial t}\|_{L^2(\Omega)}$ vs time: left
for $\Delta t=10^{-2},Re=400, T=35, dim (X_H)=6561, dim(Y_H)=1681,
dim(X_h)=25921, dim(Y_h)=6561$ and right for $\Delta
t=5\times10^{-3},Re=1000$, $T=56$ and $dim(X_H)=6561,
dim(Y_H)=1681, dim(X_h)=25921, dim(Y_h)=6561.$}
\label{FigNSstabRe4001000}
\end{figure}
\end{center}
\begin{table}[!h]
\begin{center}
\begin{tabular}{|c|c|c|c|c|}
\hline
&\multicolumn{2}{|c|}{Re=400} & \multicolumn{2}{|c|}{Re=1000}\\
\hline
\hline
Algorithm \ref{ALGO1G3} & $\Delta t=10^{-2}$ & $\Delta t=10^{-2}$ & $\Delta t=5\times10^{-3}$ & $\Delta t=5\times 10^{-3}$\\
\hline
Algorithm \ref{NSALGO1} & $\Delta t=0.1$ & $\Delta t=0.5$ & $\Delta t=0.1$ & $\Delta t=0.5$\\
  & $\tau=30$ & $\tau=30$ & $\tau=100$ & $\tau=100$ \\
\hline

$\ds\frac{Dt_{Algo \ref{NSALGO1}}}{ Dt_{III}}$ & 10 & 50 & 20 & 100 \\
\hline
%
%
\end{tabular}
\end{center}
\caption{Non incremental Navier-Stokes $\Delta t$ maximum and
stability.}\label{deltamax}
\end{table}
\noindent \newpage \noindent For small Reynolds number ($Re=100$), we do not
obtain an enhanced stability  (larger $\Delta t$) as compared to  the semi-implicit scheme. The latter
is stable even for large values of the time step $\Delta t$. When considering larger values of $Re$, says $Re=400$ and$Re=1000$,
the bi-grid scheme  appears as  10 times to
50 times more stable than the semi-implicit one.
As a consequence, the gain in CPU time to compute the steady state increases:
the stability of the semi-implicit scheme is limited by  $\Delta t=0.005$ while the bi-grid algorithm remains stable for
$\tau=100$ and $\Delta t=0.1$ and $0.5$.
\subsection{Precision test: Comparison with an exact solution }
\noindent In order to show the accuracy of the bi-grid
schemes, we simulate the Navier-Stokes equations for a Bercovier-Engelman type 
solution \cite{Bercovier}, say \[ \left.
\begin{array}{ll}
(u_1,u_2)=(-2x^2y(1-y)(1-2y)(1-x)^2\exp^{\sin
t},2y^2x(1-x)(1-2x)(1-y)^2\exp^{\sin t}) \\
p=(x-0.5)(y-0.5)
\end{array}
\right.\]
and with the corresponding term $f$. We can observe in  Figures
\ref{FigNSuexactRe1000} and \ref{FigNSerroruexactRe1000} that the solutions coincide with the reference one at final time $T=56$ and that the $L^2$ norm of the error bi-grid schemes is of the order of $\Delta t$ on all interval time.
\begin{center}
\vspace{-0.5cm}
\begin{figure}[!h]
\begin{center}
\includegraphics[height=4.6cm,width=5cm]{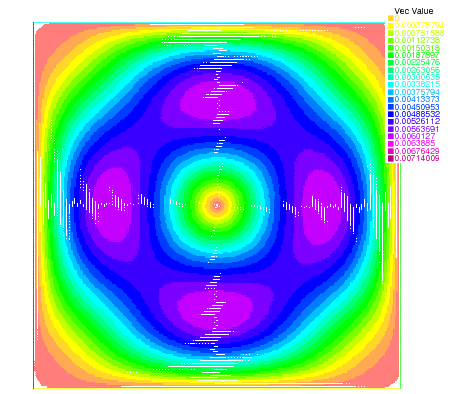}
\hspace{-0.5cm}
\includegraphics[height=4.6cm,width=5cm]{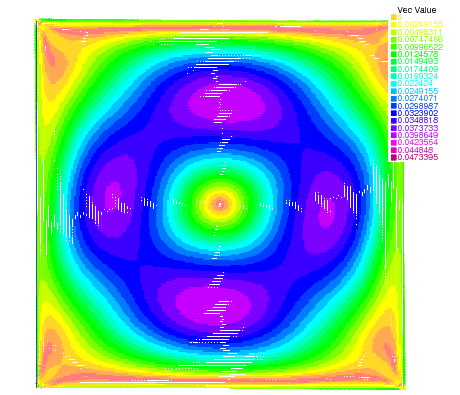}
\hspace{-0.5cm}
\includegraphics[height=4.6cm,width=5cm]{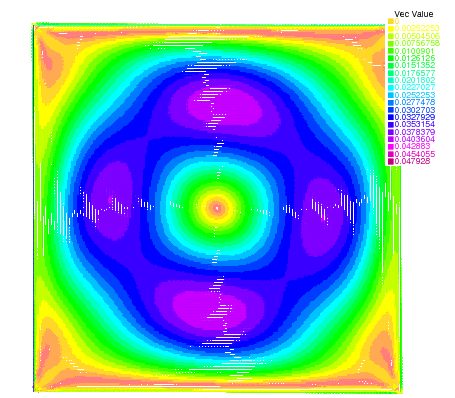}\\
%
%
\includegraphics[height=4.6cm,width=5cm]{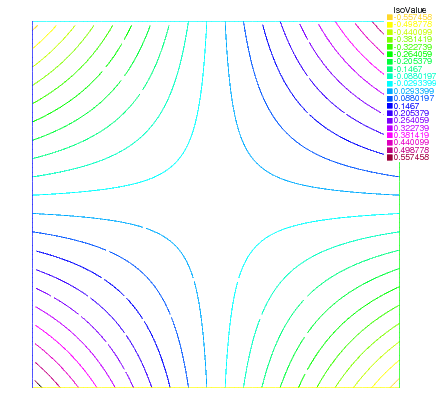}
\hspace{-0.5cm}
\includegraphics[height=4.6cm,width=5cm]{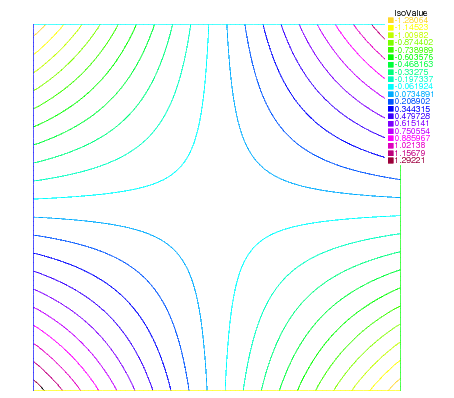}
\hspace{-0.5cm}
\includegraphics[height=4.6cm,width=5cm]{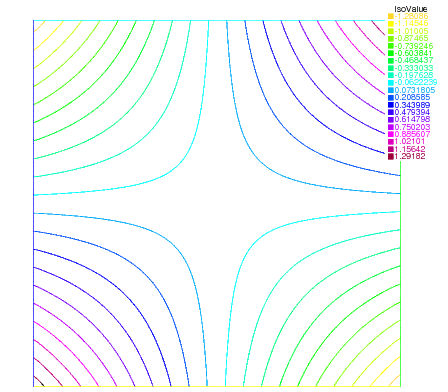}
%
\end{center}
%
\caption{From top to bottom~: velocity and pressure, and from left
to right~: the exact solution and Algorithms \ref{ALGO1G1} and
\ref{NSALGO1} solutions for $\Delta t=5\times 10^{-3},Re=1000,
T=56, dim(X_H)=6561, dim(Y_H)=1681, dim(X_h)=25921,
dim(Y_h)=6561.$} \label{FigNSuexactRe1000}
\end{figure}
\end{center}
%
\begin{center}
\begin{figure}[!h]
\vspace{-2cm}
\begin{center}
\hspace{-1.5cm}
\includegraphics[height=10.2cm,width=10.4cm]{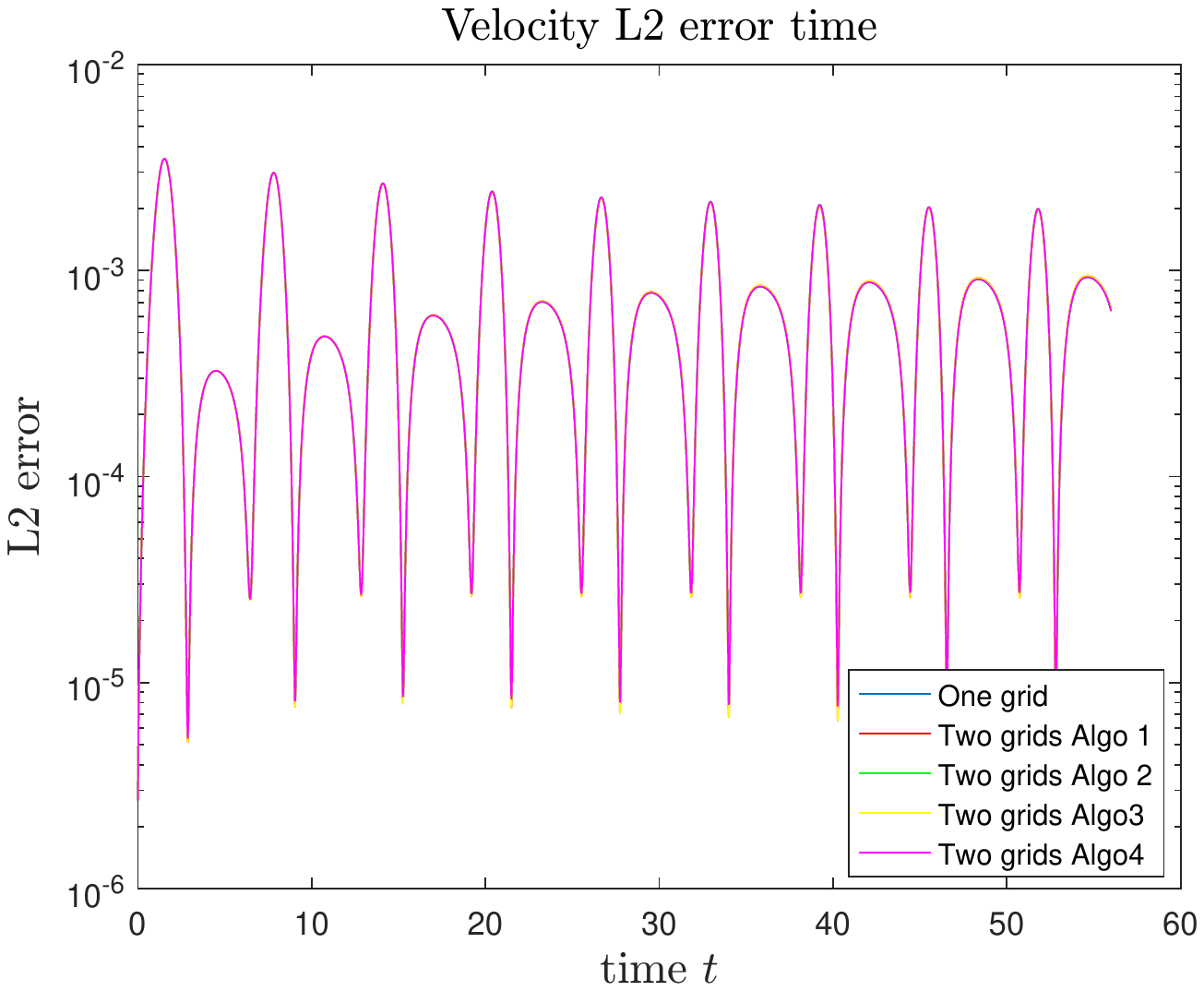}
\hspace{-3.8cm}
\includegraphics[height=10.2cm,width=10.4cm]{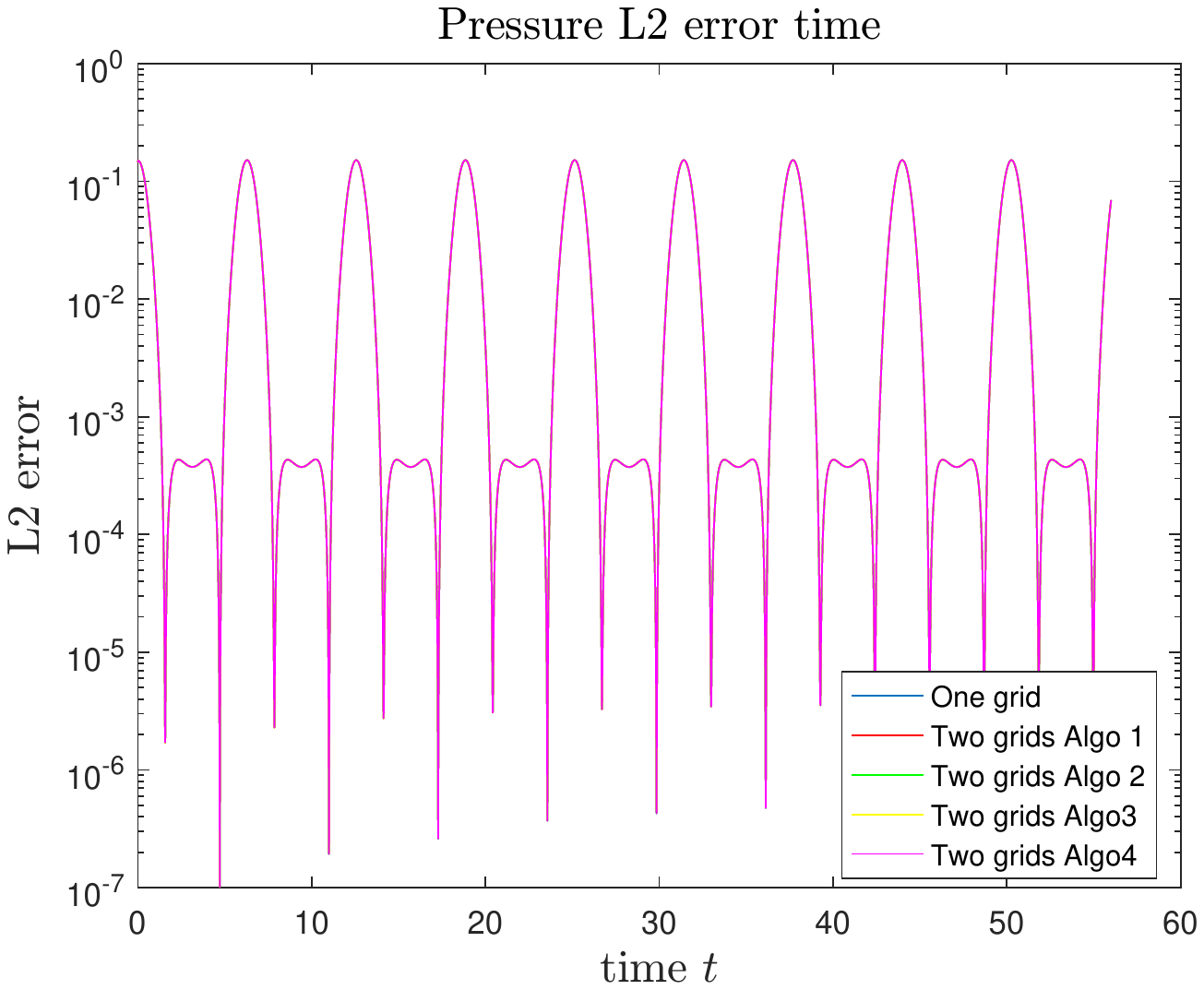}\\
\vspace{-3.3cm}
\end{center}
%
\caption{From left to right~: the $L^2-$error for the velocity and
the pressure between the exact solution and the solution of
algorithms \ref{ALGO1G1}, \ref{NSALGO1}, \ref{NSALGO2},
\ref{NSALGO3} and \ref{NSALGO4}, for $\Delta
t=5\times10^{-3},Re=1000, T=56, dim(X_H)=6561, dim(Y_H)=1681,
dim(X_h)=25921, dim(Y_h)=6561.$} \label{FigNSerroruexactRe1000}
\end{figure}
\end{center}
\newpage
\section{Concluding Remarks and perspectives}
The new bi-grid projection schemes introduced in this work allow to reduce the CPU time when compared to fully implicit projection scheme, for a comparable precision.
Their stability is enhanced as compared to semi-implicit projection scheme (larger time steps can be used). The stabilization we used on the high mode components is simple and does not deteriorate the consistency, the numerical results we obtained on the benchmark driven cavity agree with the ones of the literature, particularly the dynamics of the convergence to the steady state is close to the one of classical methods. This is a first and encouraging step before considering a strategy of multi-grid or multilevel adaptations of our methods, using more than 2 nested FEM spaces.
We have used here FEM methods for the spatial discretization, but the approach is applicable to others discretizations techniques such as finite differences or spectral methods.
%
\section*{Acknowledgment}
\noindent This project has been founded with support from the
National Council for Scientific Research in Lebanon and the
Lebanese University. We also thank the F\'ed\'eration de Recherche
ARC, CNRS FR 3399.
\section*{Conflict of interest}
\noindent  All authors declare no conflicts of interest in this paper

\end{document}